\documentclass[3p,times]{elsarticle}

\usepackage{ecrc}

\volume{}

\firstpage{1}

\journalname{}

\runauth{}

\jid{procs}

\usepackage{amsfonts,amssymb,amsmath,amsthm}
\usepackage{multicol}
\usepackage{graphicx}
\usepackage{epsfig}
\usepackage{subfigure}
\usepackage{enumerate}
\usepackage{color}
\usepackage{comment}
\usepackage{tikz} 
\usepackage[caption=false]{subfig}
\usepackage{caption}

\usepackage{lineno}

\usepackage[figuresright]{rotating}

\newtheorem{thm}{Theorem}[section] 
\newtheorem{cor}[thm]{Corollary}
\newtheorem{lem}[thm]{Lemma}
\newtheorem{prop}[thm]{Proposition}

\newtheorem{rem}{Remark}
\newtheorem{exam}[thm]{Example}

\numberwithin{equation}{section}

\begin{document}

\begin{frontmatter}

%% Title, authors and addresses

%% use the tnoteref command within \title for footnotes;
%% use the tnotetext command for the associated footnote;
%% use the fnref command within \author or \address for footnotes;
%% use the fntext command for the associated footnote;
%% use the corref command within \author for corresponding author footnotes;
%% use the cortext command for the associated footnote;
%% use the ead command for the email address,
%% and the form \ead[url] for the home page:
%%
\title{Title\tnoteref{label1}}
%% \tnotetext[label1]{}
%% \author{Name\corref{cor1}\fnref{label2}}
%% \ead{email address}
%% \ead[url]{home page}
%% \fntext[label2]{}
%% \cortext[cor1]{}
%% \address{Address\fnref{label3}}
%% \fntext[label3]{}

\dochead{}
%% Use \dochead if there is an article header, e.g. \dochead{Short communication}
%% \dochead can also be used to include a conference title, if directed by the editors
%% e.g. \dochead{17th International Conference on Dynamical Processes in Excited States of Solids}

\title{%(Compact)
Classical multivariate Hermite coordinate interpolation \textcolor{black}{on} n-dimensional grid\textcolor{black}{s}
}

%% use optional labels to link authors explicitly to addresses:
%% \author[label1,label2]{<author name>}
%% \address[label1]{<address>}
%% \address[label2]{<address>}
\author[label2]{Aristides I. Kechriniotis}
\ead[label2]{arisk7@gmail.com}
\address[label2]{Department of Physics, University of Thessaly, 3rd Km Old National Road Lamia–Athens
35100, Lamia
Greece}

\author[label1]{Konstantinos K. Delibasis\corref{cor1}}
\ead[label1]{kdelibasis@gmail.com}
\cortext[cor1]{Corresponding author. Tel.: (+30) 22310 66908.}
\address[label1]{Department of Computer Science and Biomedical Informatics, University of Thessaly, 2-4 Papasiopoulou str., P.O. 35131 Lamia, Greece}

\author[label3]{Iro Oikonomou}
\ead[label3]{iro.oikonomou99@gmail.com}

\address[label3]{Department of Informatics and Telecommunications, National and Kapodistrian University of Athens,
Panepistimioupolis, Ilisia 157 84, Athens Greece }

\author[label4]{Georgios N. Tsigaridas}
\address[label4]{Department of Physics
School of Applied Mathematical and Physical Sciences,
National Technical University of Athens,
Zografou Campus
GR-15780 Zografou, Athens
Greece}
\ead[label4]{gtsig@mail.ntua.gr}

\begin{abstract}
%% Text of abstract
In this work, we study the Hermite interpolation on $n$-dimensional non-equally spaced, rectilinear grids over a field $\Bbbk $ of characteristic zero, given the values of the function at each point of the grid and the partial derivatives up to a maximum degree. First, we prove the uniqueness of the interpolating polynomial, and we further obtain a compact closed form that uses a single summation, irrespective of the dimensionality, which is algebraically simpler than the only alternative closed form for the $n$-dimensional classical Hermite interpolation \cite{fractmult}. We provide the remainder of the interpolation in integral form; we derive the ideal of the interpolation and express the interpolation remainder using only polynomial divisions, in the case of interpolating a polynomial function. \textcolor{black}{Moreover, we prove the continuity of Hermite polynomials defined on adjacent $n-$dimensional grids, thus establishing spline behavior}. Finally, we perform illustrative numerical examples to showcase the applicability and high accuracy of the proposed interpolant, in the simple case of few points, as well as hundreds of points on 3D-grids using a spline-like interpolation, which compares favorably to state-of-the-art spline interpolation methods.

%: (The extension $\Bbbk \subset \Bbbk \left[ x_{1},x_{2},...,x_{n}% \right] /\left( f_{1}\left( x_{1}\right) ,...,f_{n}\left( x_{n}\right) \right) $ has a primitive element) $\Bbbk \left[ x_{1},x_{2},...,x_{n}\right] /\left( f_{1}\left( x_{1}\right) ,...,f_{n}\left( x_{n}\right) \right) =\Bbbk \left[ \alpha \right] $ $\ $for some $\alpha \in \Bbbk \left[x_{1},x_{2},...,x_{n}\right] /\left( f_{1}\left( x_{1}\right),...,f_{n}\left( x_{n}\right) \right) $ if and only if at most one of the univariate polynomials $f_{1},...,f_{n}$ is inseparable. This Theorem \ lead us to some results about the existence of primitive element of finite $\Bbbk -$ algebras. The above findigs are used to obtain some aswers about the following question of Frobenius: When two commuting matrices $A$ and $B~$can be written as polynomials in some matrix $C?$
\end{abstract}

\begin{keyword}
%% keywords here, in the form: keyword \sep keyword
polynomial interpolation \sep multivariate Hermite \sep classical Hermite interpolation \sep n-Dimensional grid\textcolor{black}{s}

%% PACS codes here, in the form: \PACS code \sep code
%\PACS 03.65.Pm \sep 03.50.De \sep 41.20.-q

%% MSC codes here, in the form: \MSC code \sep code
%% or \MSC[2008] code \sep code (2000 is the default)

\end{keyword}

\end{frontmatter}

%%
%% Start line numbering here if you want
%%
% \linenumbers

%% main text
%%%%%%%%%%%%%%%%%%%%%%%%%%%%%%%%%%%%%%%%%%%%%%
%Introduction
\section{Introduction and Notations}
\label{intro}
\noindent
\textcolor{black}{Interpolation is a general term for the problem of finding a function whose value at given points satisfies certain conditions. The family Hermite interpolation methods satisfy conditions involving also derivatives of the function. The method of Hermite Interpolation of the total degree is defined as using support points, considering all combinations of partial derivatives with the sum of derivative degree less than or equal to  a positive integer, for each support point separately. This interpolation is not regular (i.e. the Hermite polynomial is not unique), thus it is singular \cite{Lorenzt}.
The Hermite interpolation of the total degree is uniform if the positive integer is the same for all support points \cite[Examples 2.6 and 2.7]{sauerxu}. Hermite-Birkhoff interpolation considers different range of function derivatives for each support point, not necessarily starting from 0th order} \cite{Birkhoff}.
This work focuses on multivariate classical coordinate Hermite interpolation with support points arranged on an n-dimensional non-equally spaced rectilinear grid ($nD$ grid), given the value of a function, as well as its derivatives up to an arbitrary maximum order, defined independently for each point and each dimension. In recent years, Hermite interpolation has gained popularity among researchers; for instance, in \cite{LR,lucia2} spline interpolation was introduced with the Hermite property of having specific tangent vectors at specified points, as well as a refinement of B-splines.
% HTD is one of the least explored interpolation methods, and it utilizes the values of the function to be interpolated at a number of support points, as well as the importance of its derivatives up to an arbitrary maximum order.
An old survey of Hermite interpolation methods \cite{Lorenzt} discusses the issues of uniform, regular and singular interpolation in $n$ dimensions, however, no closed form or interpolation error is provided for the interpolant polynomial. In addition, a number of works deal with bivariate Hermite-like interpolation, where the order of partial derivatives along the two dimensions are not independent \cite{Gasca,Gevorgian}, which is different than the classical Hermite interpolation.

\textcolor{black}{In \cite{Salzerbivhyperosculatory} a special case of bivariate coordinate Hermite interpolation is discussed, called hyperosculatory interpolation over Cartesian grids, that is, interpolation problems where
all partial derivatives of first and second-order and the value of the function are known at the interpolation points 
\cite{GascaSauerhistory}.}

\textcolor{black}{While Hermite interpolation is a known subject in the mathematical community, it is rarely used in other scientific domains, including image and multidimensional signal processing. A number of tasks, such as signal resampling, image geometric transformations in 3D and frame generation in video sequences are formulated as interpolation problems and since partial derivatives up to 2nd degree are estimated via convolution with appropriate kernels, Hermite coordinate interpolation on nD-Grids is a suitable approach.
Contrary to Hermite interpolation, which requires a non-trivial generalization into multiple dimensions, most of the techniques used by the image and multidimensional signal processing community involve separable polynomial interpolation that is trivially expanded into more dimensions since the multivariate interpolating polynomial is separable. A more recent state-of-the-art such methods are collectively called generalized interpolation, which includes b-splines \cite{unser1}, \cite{unser2} and Maximal-order interpolation of minimal support -OMOMS- \cite{omom1} since the polynomials are applied not of the given values themselves, but rather on the filtered values using a cascade of appropriate causal and anti-causal infinite impulse response (IIR) filters \cite{interp_revisit} More on these methods will be discussed in the result section, where comparisons will be provided}.

We first state the generalized univariate (1D) Hermite interpolation formula proposed by
Spitzbart \cite{SPITZ}, denoting by $\mathbb{N}_{0}$ the set of non-negative integers, and by $\mathbb{N}$ the set of positive integers:
\begin{thm}\label{spy1d}
 Let $A$ be a finite subset of $\mathbb{R}$, and $\nu :A\rightarrow \mathbb{N}$
the multiplicity function. Further, let $V(A,\nu)$ be the $\mathbb{R}$-vector space $\left\{ p\in \mathbb{R} \left[ x\right] :\deg
p<\sum_{a\in A}\nu \,\left( a\right)\right\} .$ 
Given the real numbers $t_{a}^{k}, a \in A,  ~k\in \left\{
0,1,\dots,\nu \left( a\right) -1\right\} $, then there is a unique $p\in V\left( A\mathbf{,}\nu \right) $ such that $%
p^{\left( k\right) }\left( a\right) =t_{a}^{k}$, $a\in A\,,~k\in \left\{
0,1,\dots ,\nu \left( a\right) -1\right\} $, given by% 
\begin{equation}\label{polspy}
p(x)=\sum_{a\in A}\sum_{k=0}^{\nu \left( a\right) -1}t_{a}^{k}H_{a}^{k}(x), 
\end{equation}
where $\displaystyle H_{a}^{k}\left( x\right) = H_{a}\left( x\right) \frac{\left(
x-a\right) ^{k}}{k!}\sum_{t=0}^{\nu \left( a\right) -k-1}\frac{\left(
x-a\right) ^{t}}{t!}\frac{d^{t}}{dx^{t}}\left( \frac{1}{H_{a}\left( x\right) 
}\right) \left( a\right) $, with $\displaystyle H_{a}\left( x\right) :=\prod_{\substack{ %
c\in A \\ c\neq a}}\left( \frac{x-c}{a-c}\right) ^{\nu \left( a\right) }.$
\end{thm}
\noindent In the above notation, $a$ denotes any of the support points $A$ (where the values of the unknown function, as well as the values of its derivatives, are given) and the multiplicity function $v$ holds the maximum order of derivative minus 1. The above Theorem can be extended when $A$ is a subset of any field $\Bbbk $ of characteristic zero.

In our previous work \cite{DK12}, we derived a new closed-form expression of the univariate Hermite interpolating polynomial for the general case of arbitrarily spaced data, that was algebraically significantly simpler than Theorem \ref{spy1d}, since it only requires simple matrix multiplications, rather than n-order derivatives of rational polynomial functions
$H^{k}_z(x)$ of Eq.(\ref{polspy}). %$H^{k}_z(x)$ of Eq.(1). 

The 1D Hermite interpolant of Theorem \ref{spy1d} can be easily extended to multivariate interpolation. \textcolor{black}{ Namely, in \cite{GHITV} the two dimensional (2D) interpolation is considered, by constructing the fundamental Hermite polynomials using tensor products, according to Newton approach, as described in \cite{gascasa}. In \cite{DK14} we extended our previous work \cite{DK12} into two dimensions (2D), using a different generalization. Our generalization was not a Newton approach because the polynomial basis that we used was not fundamental. Consequently, the deriving closed-form expression for the polynomial was again much more compact than \cite{BGHI,GHITV}, employing matrix multiplications instead of derivatives of rational functions and applications in the case of support points arranged on a non-equidistant grid.} For both the 1D and 2D cases, we provided means for computationally efficient implementations of the proposed Hermite interpolating polynomials that achieved computational complexity comparable to other popular and much simpler interpolation techniques, such as cubic splines, whereas the measured error when applied to clinical medical images was superior.

For the generalization of Hermite interpolation into $n$ dimensions, we require the following notations.

\noindent Let $A$ be a set, $\left\vert A\right\vert $ the cardinality of $A$, and $A^{n}:=\underset{n-times}{A\times \dots\times A}$. Given the sets $A_{1},\dots,A_{n}$, then $\mathbf{A:=}A_{1}\times
\dots \times A_{n}$. Further, the element $\left( a_{1},\dots , a_{n}\right) \in 
\mathbf{A}$ will be denoted by $\mathbf{a}$. Let 
 $\mathbf{0=}\left( 0,0,\dots,0\right) $, \ $\mathbf{1=}\left(
1,1,\dots,1\right)$ be the zero vector and ones vector,  respectively.
Thus, points $\mathbf{a} = \left( a_{1},\dots, a_{n}\right)$ are arranged on a non-regular N-dimensional grid $\mathbf{A}$. Let $\Bbbk $ be a field of characteristic zero. For $\mathbf{a\in }\Bbbk ^{n}$, and $\mathbf{m\in}\mathbb{N}_{0}^{n}$ we denote $\mathbf{a}^{\mathbf{m}}:=\prod_{i=1}^{n}a_{i}^{m_{i}}$. 
Let $\mathbf{k} = \left( k_{1},\dots, k_{n}\right)$ be an $n$-dimensional vector of non-negative integers, holding the order of partial derivatives of the interpolating polynomial with respect to each variable.

In $\mathbb{N}_{0}^{n}$ we define the relation "$\leq "$ as follows: $\mathbf{k}\leq 
\mathbf{m}$ if and only if $k_{i}\leq m_{i}$, for every $i=1,\dots ,n.$ Clearly 
$\left( \mathbb{N}_{0}^{n},\leq \right) $ is a poset ($\mathbb{N}_{0}^{n}$ is partially ordered). If $\mathbf{k\leq m}$ and $\mathbf{k\neq m}$, then $\left[ \mathbf{k,m}\right] :=\left\{ \mathbf{l\in }\mathbb{N}_{0}^{n}:\mathbf{k\leq l\leq m}\right\} =\left[ k_{1},m_{1}\right] \times
\dots \times \left[ k_{n},m_{n}\right] $, and is valid $\left\vert \left[ 
\mathbf{k,m}\right] \right\vert =\prod_{i=1}^{n}\left( m_{i}-k_{i}+1\right) $.

Given the finite subsets $A_{i}$, $i=1,\dots ,n$ of the field $\Bbbk $,
and the multiplicity functions $\nu _{i}:A_{i}\rightarrow 
\mathbb{N},~i=1,\dots ,n$. Then for $i\in \left\{ 1,\dots ,n\right\} $, \textcolor{black}{ and $a_i$ any element of $A_{i}$, we define}

\textcolor{black}{
\begin{equation}\label{Hmikra}
%H_{\left( i,a\right) }\left( x_{i}\right) 
H_{a_i}\left( x_i\right) :=\prod_{\substack{ c\in A_{i}  \\ %
c\neq a_i}}\left( \frac{x_{i}-c}{a_i-c}\right) ^{\nu _{i}\left( a\right) }\in \Bbbk
\left[ x_{i}\right] .
\end{equation}
}
Let $\nu :\mathbf{A}\rightarrow 
%TCIMACRO{\U{2115} }%
%BeginExpansion
\mathbb{N}
%EndExpansion
^{n}$ be the generalized multiplicity function given by $\nu \left( 
\mathbf{a}\right) :=\left( \nu _{1}\left( a_{1}\right) ,\dots ,\nu _{n}\left(
a_{n}\right) \right) $. For $\mathbf{a\in A}$ and $\mathbf{k\in }\left[ 
\mathbf{0,}\nu \left( \mathbf{a}\right) -\mathbf{1}\right] $ we define
\begin{equation}\label{Hmegala}
H_{\left( \mathbf{a,k}\right) }\left( x_{1},\dots, x_{n}\right) :=\prod_{i=1}^{n}%
\frac{\left( x_{i}-a_{i}\right) ^{k_{i}} \textcolor{black}{H_{a_i}\left( x_i\right) }}{k_{i}!}\in \Bbbk \left[ x_{1},\dots ,x_{n}\right] .
\end{equation}

Let $V\left( \mathbf{A,}\nu \right) $ denote the $\Bbbk$-vector space $\left\{ f\in \Bbbk \left[ x_{1},\dots,x_{n}\right] :\deg
_{i}f<\sum_{a\in A_{i}}\nu _{i}\,\left( a\right) ,i=1,\dots ,n\right\} ,$ a
basis of which is the following 
\begin{equation*}
B\left( \mathbf{A,}\nu \right) :=\left\{ 
\mathbf{x}^{\mathbf{k}}:\mathbf{k\in }\left[ \mathbf{0,}-\mathbf{1+}\sum_{%
\mathbf{a\in A}}\nu \left( \mathbf{a}\right) \right] \right\} .
\end{equation*}

Further, let us define %$\displaystyle H_{a_j}\left( x_j\right) :=\prod_{\substack{ %c\in A_j \\ c\neq a_j}}\left( \frac{x_j-c}{a_j-c}\right) ^{\nu \left( a\right) }$, for $a_j \in A_j$, as well as

\begin{equation*}
    H_{a_j}^{k_j}\left( x_j\right) := H_{a_j}\left( x_j\right) \frac{\left(
x_j-a_j \right) ^{k_j}}{k_{j}!}\sum_{t=0}^{\nu \left( a_j\right) -k_j-1}\frac{\left(
x_j-a_j \right) ^{t}}{t!}\frac{d^{t}}{dx_{j}^{t}}\left( \frac{1}{H_{a_j}\left( x_j\right) 
}\right) \left( a_j\right) \in \Bbbk [x],
\end{equation*}
and subsequently, $ \displaystyle H_{\mathbf{a}}^{\mathbf{k}}\left( x_{1},\dots, x_{n}\right) :=\prod_{j=1}^{n}H_{a_j}^{k_j}(x_j) \in \Bbbk \left[
x_{1},\dots ,x_{n}\right].$  

Let us also define $P(\mathbf{A},\nu)$ as follows
\begin{equation*}
P\left( \mathbf{A},\nu \right) :=\left\{ H_{\left( \mathbf{a,k}\right) }:%
\mathbf{a\in A}\text{, }\mathbf{k\in }\left[ \mathbf{0,}\nu \left( \mathbf{a}%
\right) -\mathbf{1}\right] \right\}.
\end{equation*}

Finally, \textcolor{black}{we define the partial derivative operator acting on $f$ as}
$$
\partial ^{\mathbf{k}}:=\prod_{i=1}^{n}\partial _{i}^{k_{i}},\text{ and } \;\partial _{%
\mathbf{a}}^{\mathbf{k}}f\left( \mathbf{x}\right) :=\partial ^{\mathbf{k}%
}f\left( \mathbf{x}\right) \left\vert _{\mathbf{x=a}}\right. ,\;\text{where } \partial _{i}^{k}:=\frac{\partial ^{k}}{\partial x_{i}^{k}}.$$

To our knowledge, the only $n$-dimensional multivariate generalization of the 1D classical Hermite interpolation in Spitzbart \cite{SPITZ} was proposed very recently in \cite{fractmult}, and it is presented in the following Theorem: 
\begin{thm}\label{th1}
Let $t_{\mathbf{a}}^{\mathbf{k}}\in 
\mathbb{R}\,,~\mathbf{a\in A}:=A_{1}\times \dots \times A_{n}
,~\mathbf{k\in }\left[ \mathbf{0,}\nu \left( \mathbf{a}%
\right) -\mathbf{1}\right] .$ Then there exists a unique $p\in V\left( \mathbf{A,%
}\nu \right) $ such that $\partial _{\mathbf{a}}^{\mathbf{k}}p=t_{\mathbf{a}%
}^{\mathbf{k}},~\mathbf{a\in A},~\mathbf{k\in }\left[ \mathbf{0,}\nu \left( 
\mathbf{a}\right) -\mathbf{1}\right] $, given by 
\begin{equation}\label{ndpolmeh}
\textcolor{black}{p(x_1,\dots,x_n)=p(\mathbf{x})=}\sum_{\mathbf{a\in A}}\sum_{\mathbf{m\in }\left[ \mathbf{0,}\nu \left( 
\mathbf{a}\right) -\mathbf{1}\right] }t_{\mathbf{a}}^{\mathbf{k}}H_{\mathbf{a%
}}^{\mathbf{k}}=\sum_{a_{1}\in A_{1}} \dots \sum_{a_{n}\in
A_{n}}\sum_{k_{1}=0}^{\nu _{1}\left( a_{1}\right) -1} \dots \sum_{k_{n}=0}^{\nu
_{n}\left( a_{n}\right) -1}\ \prod_{j=1}^{n}H_{a_{j}}^{k_{j}}\left(
x_{j}\right) t_{\mathbf{a}}^{\mathbf{k}}.
\end{equation}

\begin{comment}
where 
\begin{equation}
    H_{a_j}^{k_j}\left( x_j\right) = H_{a_j}\left( x_j\right) \frac{\left(
x_j-a_j \right) ^{k_j}}{k_{j}!}\sum_{t=0}^{\nu \left( a_j\right) -k_j-1}\frac{\left(
x_j-a_j \right) ^{t}}{t!}\frac{d^{t}}{dx_{j}^{t}}\left( \frac{1}{H_{a_j}\left( x_j\right) 
}\right) \left( a_j\right) \in \Bbbk [x].
\end{equation}
with  $\displaystyle H_{a_j}\left( x_j\right) :=\prod_{\substack{ %
c\in A_j \\ c\neq a_j}}\left( \frac{x_j-c}{a_j-c}\right) ^{\nu \left( a\right) }.$
\end{comment}

\begin{comment}
\begin{equation*}
    H_{\mathbf{a%
}}^{\mathbf{k}} = \prod_{j=1}^{n}H_{a_{j}}^{k_{j}}\left(
x_{j}\right) t_{\mathbf{a}}^{\mathbf{k}}
\end{equation*}
\end{comment}
%where the symbols are defined according to the notation of  Theorem \ref{spy1d}.
\end{thm}
\textcolor{black}{
The authors of that work, also used the Newton approach to constructing the fundamental polynomial basis with tensor products of the 1D polynomial in \cite{SPITZ}, whereas the remainder of the interpolating polynomial was not provided.}
\textcolor{black}{Furthermore, the algebraic complexity of the Hermite polynomial expressed in Theorem\ref{th1} increases rapidly with the number of  dimensions $n$, which can be confirmed by comparing equations (\ref{polspy}) and (\ref{ndpolmeh})}. %Furthermore, it is easy to observe by comparing the equation \ref{spy1d} examining Eq.[--] that the algebraic complexity of the Hermite  %For instance, the expression of the 2D Hermite polynomial was significantly increased, compared to the 1D case.  
% It has to be emphasized 
%that, unlike the simple and popular interpolation techniques (e.g. spline-based), multivariate Hermite polynomials are non-separable, thus applying Hermite interpolation in higher dimensions is anything but straightforward.%

\textcolor{black}{
The present work deals with the Hermite interpolating polynomial into a regular, non-equidistant grid of support points in arbitrarily high dimensions ($nD$ grids). Our proposed generalization into n-dimensions is not based on the standard Newton approach of constructing a fundamental polynomial basis by tensor products of the polynomials in the 1D case, as described in \cite{gascasa, Sauer2006}. Rather, we define a non-fundamental polynomial basis and use reverse lexicographic ordering for the vector of derivative order to derive an elegant, compact closed-form expression with only one summation of matrix multiplications over support points and without requiring derivatives of rational polynomial functions. We also provide the interpolation remainder, which, to our knowledge, has not been given so far.}

\textcolor{black}{
Finally, we investigate the use of the proposed spline nD Hermite interpolation, regarding the continuity and the interpolation error, comparative to other state-of-the-art methods. In the majority of practical applications, such as geometric transformation of 2D and 3D images, signal resampling, frame-rate video resampling, the number of available data is so high that only spline implementation of polynomial interpolation can be applied meaningfully, to keep the total degree of the polynomial low.}

The article is organized as follows. Having defined some necessary notations, we provide in the next section the lemmas and remarks required for Theorem \ref{thm1} about the uniqueness of the interpolating polynomial and Theorem \ref{finterpform} for the proposed closed formula of the interpolant. The remainder of the interpolating polynomial is also proven in section \ref{errr}. In section \ref{idealco} the ideal of the interpolation is studied, deriving a reduced Gr{\"o}bner basis of the ideal. This basis is utilized to express the remainder of the interpolation of polynomial functions using cascaded polynomial divisions. The notations introduced in this paper, also facilitate a very short proof of Theorem \ref{th1}, which is provided in section \ref{altproof} for completeness. Finally, simple algebraic examples are given, showing that the bilinear and trilinear interpolation are special cases of the proposed n-dimensional Hermite interpolation. In addition, few arithmetic implementations of the proposed interpolation of known functions are provided, demonstrating its superiority in accuracy against other popular interpolation techniques.

\section{The proposed Hermite interpolation on $n$-$D$ rectilinear grids}
\label{mainpart}

\begin{rem}\label{rem1}
{\rm
It is easy to verify that 
\begin{eqnarray*}
\left\vert \left\{ \mathbf{a}^{\mathbf{k}}:\text{ }\mathbf{a\in A}\text{, }%
\mathbf{k\in }\left[ \mathbf{0,}\nu \left( \mathbf{a}\right) -\mathbf{1}%
\right] \right\} \right\vert =\sum_{\mathbf{a\in A}}\sum_{\mathbf{k\in }%
\left[ \mathbf{0,}\nu \left( \mathbf{a}\right) -\mathbf{1}\right] }1= 
\sum_{\mathbf{a\in A}}\prod_{i=1}^{n}\nu _{i}\left( a_{i}\right)
=\prod_{i=1}^{n}\sum_{a\in A_{i}}\nu _{i}\left( a\right) .
\end{eqnarray*}
}
\end{rem}

\begin{rem}\label{rem2}
{\rm
The cardinality of the set $P(\mathbf{A}, \nu)$ is given by $\left\vert P\left( \mathbf{A,}\nu \right) \right\vert
=\left\vert \left\{ \prod_{i=1}^{n}\left( \mathbf{x}-\mathbf{a}\right) ^{%
\mathbf{k}}:\text{ }\mathbf{a\in A}\text{, }\mathbf{k\in }\left[ \mathbf{0,}%
\nu \left( \mathbf{a}\right) -\mathbf{1}\right] \right\} \right\vert$,
which can be simplified: $$\left\vert P\left( \mathbf{A,}\nu \right) \right\vert =\sum_{%
\mathbf{a\in A}}\prod_{i=1}^{n}\nu _{i}\left( a_{i}\right).$$ The above yields $\left\vert P\left( \mathbf{A,}\nu \right) \right\vert
=\prod_{i=1}^{n}\sum_{a_{i}\in A_{i}}\nu _{i}\left( a\right) =\dim
_{k}V\left( \mathbf{A,}\nu \right)$. The set of partial derivative operators has equal cardinality: $$\left\vert \left\{ \partial _{%
\mathbf{a}}^{\mathbf{m}},~\mathbf{a\in A},~\mathbf{m\in }\left[ \mathbf{0,}%
\nu \left( \mathbf{a}\right) -\mathbf{1}\right] \right\} \right\vert
=\prod_{i=1}^{n}\sum_{a_{i}\in A_{i}}\nu _{i}\left( a\right) .$$
}
\end{rem}

\begin{rem}\label{rem3}
{\rm
Since $\deg _{i}H_{\left( \mathbf{a,k}\right) }=k_{i}+\sum_{a\in A_{i}-\left\{
a_{i}\right\} }\nu _{i}\left( a\right) <\sum_{a\in A_{i}}\nu _{i}\left(
a\right) $, the following holds $$P\left( \mathbf{A,}\nu \right) \subset V\left( 
\mathbf{A,}\nu \right).$$
}
\end{rem}
%\textcolor{black}{
\begin{rem}\label{rem4}
{\rm
For $\mathbf{a,b\in A}$, and $\mathbf{k,m\in }\left[ \mathbf{0,}\nu \left( 
\mathbf{a}\right) -\mathbf{1}\right] $ by using the Leibniz derivative rule
we easily get:
\begin{equation*}
\partial _{\mathbf{a}}^{\mathbf{k}}H_{\left( \mathbf{b,m}\right) }=\left\{ 
\begin{array}{cccccc}
0, & \text{if} & \mathbf{a\neq b} &  &  &  \\ 
0, & \text{if} & \mathbf{a=b} & \text{and} & \mathbf{k<m} &  \\ 
0, & \text{if} & \mathbf{a=b} & \text{and} & \mathbf{k,m} & \text{are
incomparable} \\ 
1, & \text{if} & \mathbf{a=b} & \text{and} & \mathbf{m=k} & 
\end{array}%
\right.
\end{equation*}
\textcolor{black}{In the case of $\mathbf{a=b}$ and $\mathbf{k>m}$, $\partial _{\mathbf{a}}^{\mathbf{k}}H_{\left( \mathbf{b,m}\right) }$ obtains non-binary value according to Eq.(\ref{Hmegala}), therefore $H_{\left( \mathbf{b,m}\right) }$ this is not a fundamental basis.}
}
\end{rem}
%}

\begin{lem}\label{poselem}
Let $\left( G,\leq \right) $ be a poset, and $A$ a finite subset of $G$. Then the system of linear equations wit unknown $\textcolor{black}{\xi}$
\begin{equation}\label{cpose}
\sum_{\beta \leq a }c_{a ,\beta }\textcolor{black}{\xi_{\beta }}=d_{a },a \in
A,  
\end{equation}%
has unique solution, where $c_{a ,\beta }\in \Bbbk $ with $c_{a
,a }\neq 0$ for all $a \in A$.
\end{lem}

\begin{proof}
Let $\left\vert A\right\vert =m$. Then we number the elements of $A$ as
follows: We choose radomly a minimal element of $A\,$which will be denoted
by $a_{1}.$ Next we choose a minimal element of $A-\left\{ a
_{1}\right\} $ denoted by $a_{2}$, and so on. That is $A=\left\{
a_{1},\dots,a_{m}\right\} .$ Therefore (\ref{cpose}) can
be rewritten as 
\begin{equation*}
\sum_{i=1}^{j}\varepsilon _{i,j}c_{a_i}\textcolor{black}{\xi_{a_{i}}}=d_{a
_{j}},j=1,\dots,m,
\end{equation*}%
where 
\begin{equation*}
   \varepsilon _{i,j}:= \begin{cases}
     0, & \; \text{if }\; a _{j}<a_{i}  \\
     0, & \; \text{if }\; a _{i},a_{j} \text{ are not comparable}\\
     1, & \; \text{if }\;  a_{i}\leq a_{j},
   \end{cases}
\end{equation*}

Clearly, the matrix of the coefficients of the variables $x_{a_{i}}$
is lower triagular, and because all its diagonal elements $\varepsilon
_{j,j}c_{a_{j}}$ are non zero, we conclude that (\ref{cpose})
has exactly one solution.
\end{proof}

We now present the Theorem that  \textcolor{black}{Hermite interpolating polynomial of coordinate degree on multivariate grids is regular}.
\begin{thm}\label{thm1}
Given the elements $t_{\mathbf{a}}^{\mathbf{k}}\in \Bbbk \,,~\mathbf{a\in A}%
,~\mathbf{k\in }\left[ \mathbf{0,}\nu \left( \mathbf{a}\right) -\mathbf{1}%
\right] $, there exists a unique $f\in V\left( \mathbf{A,}\nu \right) $
such that $\partial _{\mathbf{a}}^{\mathbf{m}}f=t_{\mathbf{a}}^{\mathbf{m}},~%
\mathbf{a\in A},~\mathbf{m\in }\left[ \mathbf{0,}\nu \left( \mathbf{a}%
\right) -\mathbf{1}\right]$.
\end{thm}

\begin{proof}
It is sufficient to show that there are unique elements $\textcolor{black}{\xi_{\mathbf{a}}^{\mathbf{%
m}}}\in k,~\mathbf{a\in A,~m\in }\left[ \mathbf{0,}\nu \left( \mathbf{a}%
\right) -\mathbf{1}\right] $ , such that the polynomial 
\begin{equation}\label{frel1}
f(\mathbf{x})=\sum_{\mathbf{b\in A}}\sum_{\mathbf{k\in }\left[ \mathbf{0,}\nu \left( 
\mathbf{b}\right) -\mathbf{1}\right] }\textcolor{black}{\xi_{\mathbf{b}}^{\mathbf{k}}}H_{\left( 
\mathbf{b,k}\right) } (\mathbf{x})
\end{equation}%
satisfies the conditions 
\begin{equation*}
\partial _{\mathbf{a}}^{\mathbf{m}}f=t_{\mathbf{a}}^{\mathbf{m}}.
%,~\mathbf{m\in }\left[ \mathbf{0,}\nu \left( \mathbf{a}\right) -\mathbf{1}\right].
\end{equation*}%
By applying the derivative operators $\partial _{\mathbf{a}}^{\mathbf{m}}$ , $%
\mathbf{a\in A,~m\in }\left[ \mathbf{0,}\nu \left( \mathbf{a}\right) -%
\mathbf{1}\right] $ to  (\ref{frel1}), we get the following system
consisting of a number of $ \sum_{\mathbf{a\in A}}\prod_{i=1}^{n}\nu _{i}\left(
a_{i}\right) $  linear equations with an equal number of variables $\textcolor{black}{\xi_{\mathbf{a}%
}^{\mathbf{k}}},$ $\mathbf{a\in A,~m\in }\left[ \mathbf{0,}\nu \left( \mathbf{%
a}\right) -\mathbf{1}\right] $: 
\begin{equation*}
\partial _{\mathbf{a}}^{\mathbf{m}}f=\sum_{\mathbf{b\in A}}\sum_{\mathbf{%
k\in }\left[ \mathbf{0,}\nu \left( \mathbf{b}\right) -\mathbf{1}\right] }\textcolor{black}{\xi_{%
\mathbf{b}}^{\mathbf{k}}}\partial _{\mathbf{a}}^{\mathbf{m}}H_{\left( \mathbf{%
b,k}\right) },
\end{equation*}%
or 
\begin{equation*}
\sum_{\mathbf{b\in A}}\sum_{\mathbf{k\in }\left[ \mathbf{0,}\nu \left( 
\mathbf{b}\right) -\mathbf{1}\right] }\textcolor{black}{\textcolor{black}{\xi_{\mathbf{b}}^{\mathbf{k}}}}\partial _{%
\mathbf{a}}^{\mathbf{m}}H_{\left( \mathbf{b,k}\right) }=t_{\mathbf{a}}^{%
\mathbf{m}},
\end{equation*}%
which by Remark \ref{rem3} can be partitioned to the following linear equation
systems %$\left( S_{\mathbf{a}}\right) ,$ 
\begin{equation}\label{sysrel1}
\sum_{\mathbf{k\in }\left[ \mathbf{0,m}\right] }\textcolor{black}{\xi_{\mathbf{a}}^{\mathbf{k}}%
}\partial _{\mathbf{a}}^{\mathbf{m}}H_{\left( \mathbf{a,k}\right) }=t_{%
\mathbf{a}}^{\mathbf{m}},\; \mathbf{m\in }\left[ \mathbf{0,}\nu \left( \mathbf{%
a}\right) -\mathbf{1}\right].  
\end{equation}%
Clearly each one of the linear systems (\ref{sysrel1}) has a number of $\left\vert \nu \left( \mathbf{a}%
\right) \right\vert =\prod_{i=1}^{n}\nu _{i}\left( a_{i}\right) $ variables, $\textcolor{black}{\xi_{\mathbf{a}}^{\mathbf{k}}}$, and consists of $\left\vert \nu
\left( \mathbf{a}\right) \right\vert $ equations. Since $\left[ \mathbf{0,}%
\nu \left( \mathbf{a}\right) -\mathbf{1}\right] $ is a finite subset of the
partially ordered $%
%TCIMACRO{\U{2124} }%
%BeginExpansion
\mathbb{Z}
%EndExpansion
^{n}$, and because, by Remark \ref{rem4} for each $\mathbf{m\in }\left[ \mathbf{0,}%
\nu \left( \mathbf{a}\right) -\mathbf{1}\right] $ is valid $\partial _{%
\mathbf{a}}^{\mathbf{m}}H_{\left( \mathbf{a,m}\right) }=1,$then by Lemma \ref{poselem}
follows that the system (\ref{frel1}) has exact a solution.
\end{proof}

As a result of Theorem \ref{thm1} we immediately obtain the following Corollaries, which will be used in Section 4.

\begin{cor}
The set $P\left( \mathbf{A,}\nu \right) $ is a basis of the $\Bbbk$-vector
space $V\left( \mathbf{A,}\nu \right) .\,$
\end{cor}

\begin{cor}\label{cor24}
If $f$ is an element of $V\left( \mathbf{A,}\nu \right) $, such that $%
\partial _{\mathbf{a}}^{\mathbf{m}}f=0$ for each$~\mathbf{a\in A}$,$~\mathbf{%
m\in }\left[ \mathbf{0,}\nu \left( \mathbf{a}\right) -\mathbf{1}\right] $, then $f=0.$
\end{cor}

\begin{comment}
\textcolor{black}{
\begin{cor}
The $\Bbbk$-algebra $k\left[ x_{1},\dots ,x_{n}\right] /J\left( \mathbf{A,}\nu
\right) $, as a $\Bbbk$-vector space, is isomorphic to $V\left( \mathbf{A,}%
\nu \right)$ and the following holds:
\begin{equation*}
\dim k\left[ x_{1},\dots,x_{n}\right] /J\left( \mathbf{A,}\nu \right)
=\prod_{i=1}^{n}\sum_{a_{i}\in A_{i}}\nu _{i}\left( a\right) 
\end{equation*}
\end{cor}
}

\begin{proof}
Considering Remarks \ref{rem1} and \ref{rem2}, it is sufficient to prove that the elements of $%
P\left( \mathbf{A,}\nu \right) $ are linearly independent. Let $\sum_{\mathbf{%
b\in A}}\sum_{\mathbf{k\in }\left[ \mathbf{0,}\nu \left( \mathbf{b}\right) -%
\mathbf{1}\right] }x_{\mathbf{b}}^{\mathbf{k}}H_{\left( \mathbf{b,k}\right)
}=0$, where $x_{\mathbf{b}}^{\mathbf{k}}\in \Bbbk .$ Then for any $\mathbf{%
a\in A,}$ similarly to Theorem \ref{thm1}, we obtain the following homogeneous
system with variables $x_{\mathbf{b}}^{\mathbf{k}}$:
\begin{equation*}
\sum_{\mathbf{b\in A}%
}\sum_{\mathbf{k\in }\left[ \mathbf{0,}\nu \left( \mathbf{b}\right) -\mathbf{%
1}\right] }x_{\mathbf{b}}^{\mathbf{k}}\partial _{\mathbf{a}}^{\mathbf{m}%
}H_{\left( \mathbf{b,k}\right) }=0,\mathbf{~m\in }\left[ \mathbf{0,}\nu
\left( \mathbf{a}\right) -\mathbf{1}\right], 
\end{equation*}
which has exactly one
solution: $x_{\mathbf{a}}^{\mathbf{m}}=0,~\mathbf{m\in }\left[ \mathbf{0,}%
\nu \left( \mathbf{a}\right) -\mathbf{1}\right] .$
\end{proof}

\end{comment}

We will derive an expression for the interpolating polynomial $\ f$
in Theorem \ref{thm1}. First we will use the degree reverse lexicographic order,
which will be denoted by $\prec $ . More specifically,$(k_1,\dots, k_n) \prec (l_1,\dots l_n)$, if either of the following holds:
\begin{itemize}
    \item[(i)] $k_1 +\dots + k_n < l_1 +\dots +l_n$, or
    \item[(ii)] $k_1 +\dots +k_n = l_1 +\dots+l_n$ and $k_i > l_i$ for the largest $i$ for which $k_i \neq l_i$.
\end{itemize}
For example, the reverse lexicographic order  of the elements in $\left[ \mathbf{0,}\nu \left( \mathbf{a}\right) -\mathbf{1}\right] $ is
\begin{eqnarray*}
1_{\mathbf{a}} :&=&\left( 0,\dots,0\right) \prec 2_{\mathbf{a}}:=\left(
0,\dots,0,1\right) \prec 3_{\mathbf{a}}:=\left( 0,\dots,0,1,0\right) \prec
\dots\prec \left( n+1\right) _{\mathbf{a}}:=\left( 1,0,\dots,0\right) \\
&\prec &\left( n+2\right) _{\mathbf{a}}:=\left( 0,\dots,0,2\right) \prec
\dots\prec ({\left\vert \nu \left( \mathbf{a}\right) \right\vert -1})_{\mathbf{a}}:=\left(
\nu _{1}\left( a_{1}\right) -1,\dots,\nu _{n}\left( a_{n}\right) -1\right) .
\end{eqnarray*}%
Note that from $\mathbf{m\leq n}$ follows $\mathbf{m\preceq n}$, and from $%
\mathbf{m<n}$ follows $\mathbf{m\prec n}$. That means $\mathbf{\preceq }$ is
a linear extension of $\mathbf{\leq }$.

The following Theorem provides the closed form of the interpolating polynomial, \textcolor{black} {which is not a Newton approach \cite[subsection 2.6]{Davis1975}, as it stems from Remark \ref{rem4}, although it may resemble it}.
\begin{thm}\label{finterpform}
The formula of the interpolating
polynomial $\ f$ in (\ref{frel1}) is the following
\begin{equation} \label{mainform}
f(\mathbf{x}):=\sum\limits_{\mathbf{a\in A}}\Lambda _{\mathbf{a}}^{-1}T_{\mathbf{a}}H_{%
\mathbf{a}}=
\textcolor{black}{\sum\limits_{\mathbf{a\in A}}\sum_{i=1}^{\left\vert \nu \left( 
\mathbf{a}\right) \right\vert -1 }\left( I_{\left\vert \nu \left( \mathbf{a}%
\right) \right\vert }-\Lambda _{\mathbf{a}}\right)^i T_{\mathbf{a}}H_{\mathbf{a%
}}},
\end{equation}%
\textcolor{black}{where $H_{\mathbf{a}}=H_{\mathbf{a}}(\mathbf{x})$},
\begin{eqnarray}\label{La}
\Lambda _{\mathbf{a}} &=&\left[ 
\begin{array}{ccccc}
1 & 0 & \cdots &  0 & 0 \\ 
%\varepsilon _{2,1}
\varepsilon _{2,1}\partial _{\mathbf{a}}^{2_{\mathbf{a}}}H_{\left( \mathbf{a,%
}1_{\mathbf{a}}\right) } & 1 & \cdots & 0 & 0 \\ 
\vdots & \vdots & \ddots  & \vdots & \vdots \\ 
%\varepsilon _{\left\vert \nu \left( \mathbf{a}\right) \right\vert-1,1}
\varepsilon _{\left\vert \nu \left( \mathbf{a}\right) \right\vert-1
,1} \partial _{\mathbf{a}}^{\left( \left\vert \nu \left( \mathbf{a}\right)
\right\vert -1\right) _{\mathbf{a}}}H_{\left( \mathbf{a,}1_{\mathbf{a}%
}\right) } 
& 
\varepsilon_{\left\vert \nu \left( \mathbf{a}\right)\right\vert -1,2}
\partial _{\mathbf{a}}^{\left( \left\vert \nu \left( 
\mathbf{a}\right) \right\vert -1\right) _{\mathbf{a}}}H_{\left( \mathbf{a,}%
2_{\mathbf{a}}\right) } 
&\cdots & 1 & 0 \\ 
\varepsilon _{\left\vert \nu \left( \mathbf{a}\right) \right\vert,1}
\partial _{\mathbf{a}}^{\left\vert \nu \left( \mathbf{a}\right)
\right\vert _{\mathbf{a}}}H_{\left( \mathbf{a,}1_{\mathbf{a}}\right) } 
& 
\varepsilon _{\left\vert \nu \left( \mathbf{a}\right) \right\vert,2}
\partial _{\mathbf{a}}^{\left\vert \nu \left( \mathbf{a}\right)
\right\vert _{\mathbf{a}}}H_{\left( \mathbf{a,}2_{\mathbf{a}}\right) } & 
\cdots & 
\varepsilon _{\left\vert \nu \left( \mathbf{a}\right) \right\vert,\left\vert \nu \left( \mathbf{a}\right) \right\vert -1}
\partial _{\mathbf{a}%
}^{\left\vert \nu \left( \mathbf{a}\right) \right\vert _{\mathbf{a}%
}}H_{\left( \mathbf{a,}\left( \left\vert \nu \left( \mathbf{a}\right)
\right\vert -1\right) _{\mathbf{a}}\right) } & 1%
\end{array}%
\right] , \\
x_{\mathbf{a}} &=&\left[ 
\begin{array}{c}
x_{\mathbf{a}}^{1_{\mathbf{a}}} \\ 
\vdots
\\
x_{\mathbf{a}}^{\left\vert \nu \left( \mathbf{a}\right) \right\vert }%
\end{array}%
\right] \text{, }T_{\mathbf{a}}=\left[ 
\begin{array}{c}
t_{\mathbf{a}}^{1_{\mathbf{a}}} \\ 
\vdots
\\
t_{\mathbf{a}}^{\left\vert \nu \left( \mathbf{a}\right) \right\vert }%
\end{array}%
\right] ,H_{\mathbf{a}}=\left[ 
\begin{array}{c}
H_{\left( \mathbf{a,}1_{\mathbf{a}}\right) } \\ 
\vdots
\\
H_{\left( _{\mathbf{a,}\left\vert \nu \left( \mathbf{a}\right) \right\vert
}\right) }%
\end{array}%
\right] 
,\; \lvert \nu \left( \mathbf{a}\right) \rvert =\prod_{i=1}^{n}\nu _{i}\left( a_{i}\right).\nonumber 
\end{eqnarray}
\end{thm}

\begin{proof}
The interpolating
polynomial $\ f$ in (\ref{frel1}) can be equivalently written as%
\begin{equation*}
f=\textcolor{black}{\sum\limits_{\mathbf{a\in A}}\sum_{i=0}^{\left\vert \nu \left( \mathbf{a}%
\right) \right\vert-1 }}\textcolor{black}{\textcolor{black}{\xi_{\mathbf{a}}^{i_{\mathbf{a}}}}}H_{\left( \mathbf{a,}i_{%
\mathbf{a}}\right) }
\end{equation*}%
or%
\begin{equation}\label{frel2}
f=\sum\limits_{\mathbf{a\in A}}\textcolor{black}{\textcolor{black}{\xi_{\mathbf{a}}}H_{\mathbf{a}}},  
\end{equation}%
where $\textcolor{black}{\textcolor{black}{\xi_{\mathbf{a}}}}=\left[ \textcolor{black}{\textcolor{black}{\xi_{\mathbf{a}}^{i_{\mathbf{a}}}}}\right] $ is a
column matrix of size $\left\vert \nu \left( \mathbf{a}\right) \right\vert
\times 1$, and $H_{\mathbf{a}}=\left[ H_{\left( \mathbf{a,}i_{\mathbf{a}%
}\right) }\right] $ is a row matrix of size $1\times \left\vert \nu \left( 
\mathbf{a}\right) \right\vert .$ Then the system of linear equations
becomes, 

\begin{equation*}
\sum_{i=1}^{j}\varepsilon _{i,j}\partial _{\mathbf{a}}^{j_{\mathbf{a}%
}}H_{\left( \mathbf{a,}i_{\mathbf{a}}\right) }\textcolor{black}{\xi_{\mathbf{a}}^{i_{\mathbf{a}}%
}}=t_{\mathbf{a}}^{j\mathbf{_{\mathbf{a}}}},\mathbf{~}j=1,\dots ,\left\vert \nu
\left( \mathbf{a}\right) \right\vert ,
\end{equation*}%
where 
\begin{equation*}
   \varepsilon _{i,j}:= \begin{cases}
     0, & \; \text{if }\; j_{\mathbf{a}}<i_{\mathbf{a}}  \\
     0, & \; \text{if }\; i_{\mathbf{a}},\ j_{\mathbf{a}} \text{ are not comparable}\\
     1, & \; \text{if }\;  i_{\mathbf{a}}\leq j_{\mathbf{a}},
   \end{cases}
\end{equation*}

Now we write this system of linear equations in matrix form: 
\begin{equation}\label{larel1}
\Lambda _{\mathbf{a}}\textcolor{black}{\xi_{\mathbf{a}}}=T_{\mathbf{a}},  
\end{equation}%
where $\Lambda _{\mathbf{a}}=\left[ \varepsilon _{i,j}\partial _{\mathbf{a}%
}^{j_{\mathbf{a}}}H_{\left( \mathbf{a,}i_{\mathbf{a}}\right) }\right] $ is a
lower unitriagular matrix of size $\left\vert \nu \left( \mathbf{a}\right)
\right\vert \times \left\vert \nu \left( \mathbf{a}\right) \right\vert $,
and $T_{\mathbf{a}}=\left[ t_{\mathbf{a}}^{i_{\mathbf{a}}}\right] $ is a
column matrix of size $\left\vert \nu \left( \mathbf{a}\right) \right\vert
\times 1$. \ From (\ref{larel1}) we have $\textcolor{black}{\xi_{\mathbf{a}}}=\Lambda _{%
\mathbf{a}}^{-1}T_{\mathbf{a}}$, or by Lemma 2.1 in \textcolor{black}{ref \cite{DK12}}, 
\begin{equation}\label{xarel1}
\textcolor{black}{\xi_{\mathbf{a}}=\sum_{i=0}^{\left\vert \nu \left( \mathbf{a}\right)
\right\vert -1}\left( I_{\left\vert \nu \left( \mathbf{a}\right) \right\vert
}-\Lambda _{\mathbf{a}}\right) ^{i}T_{\mathbf{a}}},  
\end{equation}%
where $I_{\left\vert \nu \left( \mathbf{a}\right) \right\vert }$ is unit
matrix of size $\left\vert \nu \left( \mathbf{a}\right) \right\vert $.
Finally, substituting (\ref{xarel1}) in (\ref{frel2}) yileds the required formula.
\end{proof}

%\section{Error of interpolation}
\section{Remainder of the interpolation for $\Bbbk=\mathbb{R}$ or $\mathbb{C}$}
\label{errr}
In this section, we derive an expression for the error of the interpolation formula given
in Theorem \ref{finterpform}, which is identical to the error of Theorem \ref{th1}, for the class of real functions   $f\left( x_{1},\dots,x_{n}\right) $, which can be continued  analytically as a
single valued, \textcolor{black}{regular functions $f\left( z_{1},\dots,z_{n}\right) $} of $n$
complex variables in a certain cross-product region $D_{z_{1}}\times
\dots \times D_{z_{n}}$.

In the case of a single variable $x$, let $C$ be a closed contour
in the region $D_{z}$ of analytic continuation of a real function $f(x)$
containing the points $a\in A$ in its interior. Let as denote 
\begin{equation*}
H\left( x\right) =\prod_{a\in A}\left( x-{a}\right) ^{\nu \left( a\right) }.
\end{equation*}

By applying the residue theorem to the contour integral  
\begin{equation*}
\frac{1}{2\pi i}\int_{C}\frac{f\left( z\right) }{\left( z-x\right) H\left(
z\right) }dz,
\end{equation*}%
we obtain (see \cite{GHITV}) 
\begin{equation}\label{fpxint1}
f(x)=p\left( x\right) +\frac{H\left( x\right) }{2\pi i}\int_{C}\frac{f\left(
z\right) }{\left( z-x\right) H\left( z\right) }dz,
\end{equation}%
where 
\begin{equation*}
p\left( x\right) =\sum_{a\in A}\sum_{k=0}^{\nu \left( a\right) -1}f^{\left(
k\right) }\left( a\right) H_{a}^{k},
\end{equation*}

From (\ref{fpxint1}) we have 
\begin{equation*}
\frac{1}{2\pi i}\int_{C}\frac{f\left( z\right) }{z-x}dz=p\left( x\right) -%
\frac{H\left( x\right) }{2\pi i}\int_{C}\frac{f\left( z\right) }{\left(
z-x\right) H\left( z\right) }dz,
\end{equation*}%
or equivalently 
\begin{equation}\label{fp1}
\frac{1}{2\pi i}\int_{C}\frac{H\left( z\right) -H\left( x\right) }{H\left(
z\right) \left( z-x\right) }f\left( z\right) dz=p\left( x\right).
\end{equation}

In the following, we assume that $C_{j},\;j=1,\dots,n$ are simple closed contours in the
regions $D_{j}$  containing $A_{j}$ in their
interior, of analyticity of $f\left(x_{1},\dots,x_{j-1},z_{j},x_{j+1},\dots,x_{n}\right) $, where $x_{1},\dots,x_{j-1},x_{j+1},\dots,x_{n}$ are fixed. Further, we assume that $f\left( z_{1},\dots,z_{n}\right) $ is \textcolor{black}{holomorphic \cite{Vladimirov2007,Fisher1999} , i.e.} 
simultaneously analytic in $D_{z_{1}}\times \dots\times D_{z_{n}}$.

Finally, we introduce the following notation: 
\begin{equation*}
H_r (w_r) = \prod_{a\in A_r} (w_r -a)^{\nu_r (a)},\; r=1,\dots, n.
\end{equation*}

\begin{lem}\label{lemErr}
The following identity holds:%
\begin{eqnarray}\label{lemrel1}
&&\frac{1}{\left( 2\pi i\right) ^{n}}\int_{C_{1}}\dots\int_{C_{n}}\frac{%
\prod_{r=1}^{n}\left( H_{r}\left( z_{r}\right) -H_{r}\left( x_{r}\right)
\right) }{\prod_{r=1}^{n}\left( z_{r}-x_{r}\right) H_{r}\left( z_{r}\right) }%
f\left( z_{1},\dots,z_{n}\right) dz_{1}\dots dz_{n} =\sum_{\mathbf{a\in A}}\sum_{\mathbf{k\in }\left[ \mathbf{0,}\nu \left( 
\mathbf{a}\right) -\mathbf{1}\right] }H_{\mathbf{a}}^{\mathbf{k}}\partial _{%
\mathbf{a}}^{\mathbf{k}}f.  
\end{eqnarray}
\end{lem}

\begin{proof}
We will prove by induction with respect to $n$: For $n=1$, (\ref{lemrel1}) reduces to (\ref{fpxint1}). Let (\ref{lemrel1}) hold for $n-1$. Then 
\begin{eqnarray*}
&&\frac{1}{\left( 2\pi i\right) ^{n}}\int_{C_{1}}\dots\int_{C_{n}}\frac{%
\prod_{r=1}^{n}\left( H_{r}\left( z_{r}\right) -H_{r}\left( x_{r}\right)
\right) }{\prod_{r=1}^{n}\left( z_{r}-x_{r}\right) H_{r}\left( z_{r}\right) }%
f\left( z_{1},\dots,z_{n}\right) dz_{1}\dots dz_{n} \\
&=&\frac{1}{2\pi i}\int_{C_{n}}\frac{H_{n}\left( z_{n}\right) -H_{n}\left(
x_{n}\right) }{\left( z_{n}-x_{n}\right) H_{n}\left( z_{r}\right) }\left( 
\frac{1}{\left( 2\pi i\right) ^{n-1}}\int_{C_{1}}\dots\int_{C_{n-1}}\frac{%
\prod_{r=1}^{n-1}\left( H_{r}\left( z_{r}\right) -H_{r}\left( x_{r}\right)
\right) }{\prod_{r=1}^{n-1}\left( z_{r}-x_{r}\right) H_{r}\left(
z_{r}\right) }f\left( z_{1},\dots,z_{n}\right) dz_{1}\dots dz_{n-1}\right) dz_{n}
\\
&=&\frac{1}{2\pi i}\int_{C_{n}}\frac{H_{n}\left( z_{n}\right) -H_{n}\left(
x_{n}\right) }{\left( z_{n}-x_{n}\right) H_{n}\left( z_{r}\right) }%
\sum_{a_{n}\in A_{1}}\dots\sum_{a_{n}\in A_{n}}\sum_{k_{1}=0}^{\nu _{1}\left(
a_{1}\right) -1}\dots\sum_{k_{n-1}=0}^{\nu _{n-1}\left( a_{n-1}\right)
-1}\prod_{j=1}^{n-1}H_{a_{j}}^{k_{j}}\left( x_{j}\right) \frac{\partial ^{%
%\NEG
{k}_{1}+\dots+k_{n-1}}}{\partial x_{1}^{k_{1}}\dots \partial x_{n-1}^{k_{n-1}}%
}f\left( a_{1,}\dots, a_{n-1},z_{n}\right) dz_{n}\\
&=& \sum_{a_{n}\in A_{1}}\dots\sum_{a_{n}\in A_{n}}\sum_{k_{1}=0}^{\nu _{1}\left(
a_{1}\right) -1}\dots\sum_{k_{n-1}=0}^{\nu _{n-1}\left( a_{n-1}\right)
-1}\prod_{j=1}^{n-1} H_{a_{j}}^{k_{j}}(x_j)\left( \frac{1}{2\pi i} \int_{C_{n}}\frac{H_{n}\left( z_{n}\right) -H_{n}\left(
x_{n}\right) }{\left( z_{n}-x_{n}\right) H_{n}\left( z_{r}\right) }  \frac{\partial ^{%
%\NEG
{k}_{1}+\dots+k_{n-1}}}{\partial x_{1}^{k_{1}}\dots \partial x_{n-1}^{k_{n-1}}%
}f\left( a_{1,}\dots, a_{n-1},z_{n}\right) dz_{n}  \right) 
.
\end{eqnarray*}%
Thus, using (\ref{fp1}) we obtain:  
\begin{eqnarray*}
&&\frac{1}{\left( 2\pi i\right) ^{n}}\int_{C_{1}}\dots\int_{C_{n}}\frac{%
\prod_{r=1}^{n}\left( H_{r}\left( z_{r}\right) -H_{r}\left( x_{r}\right)
\right) }{\prod_{r=1}^{n}\left( z_{r}-x_{r}\right) H_{r}\left( z_{r}\right) }%
f\left( z_{1},\dots,z_{n}\right) dz_{1}\dots dz_{n} \\
&=&\sum_{a_{1}\in A_{1}}\dots\sum_{a_{n}\in A_{n}}\sum_{k_{n}=0}^{\nu
_{1}\left( a_{1}\right) -1}\dots \sum_{k_{n}=0}^{\nu _{n}\left( a_{n}\right)
-1}\ \prod_{j=1}^{n}H_{a_{j}}^{k_{j}}\left( x_{j}\right) \frac{\partial ^{%
%\NEG
{k}_{n}+\dots+k_{1}}}{\partial x_{n}^{k_{n}}\dots\partial x_{1}^{k_{1}}}%
f\left( a_{1},\dots, a_{n}\right) \\
&=&\sum_{\mathbf{a\in A}}\sum_{\mathbf{k\in }\left[ \mathbf{0,}\nu \left( 
\mathbf{a}\right) -\mathbf{1}\right] }H_{\mathbf{a}}^{\mathbf{k}}\partial _{%
\mathbf{a}}^{\mathbf{k}}f.
\end{eqnarray*}
\end{proof}

\textcolor{black}{The above Lemma is a generalization into n-dimensions of Corollary 3.6.2 in \cite{Davis1975}. The interpolation error is derived in the following Theorem, which is also an n-dimensional generalization of Theorem 3.6.1 of \cite{Davis1975}}.
\begin{thm}\label{errform}
The error of the interpolation formula of Theorem \ref{finterpform} is defined as
\begin{equation*}
R:=f-\sum_{\mathbf{a\in A}}\sum_{\mathbf{k\in }\left[ \mathbf{0,}\nu \left( 
\mathbf{a}\right) -\mathbf{1}\right] }t_{\mathbf{a}}^{\mathbf{k}}H_{\mathbf{a%
}}^{k},
\end{equation*}%
where $t_{\mathbf{a}}^{\mathbf{k}}=\partial _{\mathbf{a}}^{\mathbf{k}}f$ ,
 is given by:
\begin{equation}\label{211}
R=\frac{1}{\left( 2\pi i\right) ^{n}}\int_{C_{1}}\dots \int_{C_{n}}\frac{%
\prod_{r=1}^{n}H_{r}\left( z_{r}\right) -\prod_{r=1}^{n}\left( H_{r}\left(
z_{r}\right) -H_{r}\left( x_{r}\right) \right) }{\prod_{r=1}^{n}\left(
z_{r}-x_{r}\right) H_{r}\left( z_{r}\right) }f\left( z_{1},\dots,z_{n}\right)
dz_{1}\dots dz_{n}. 
\end{equation}
\end{thm}

\begin{proof}
From (\ref{lemrel1}) we have
\begin{eqnarray}\label{212}
&&\frac{1}{\left( 2\pi i\right) ^{n}}\int_{C_{1}}\dots\int_{C_{n}}\frac{%
\prod_{r=1}^{n}H_{r}\left( z_{r}\right) -\prod_{r=1}^{n}\left( H_{r}\left(
z_{r}\right) -H_{r}\left( x_{r}\right) \right) }{\prod_{r=1}^{n}\left(
z_{r}-x_{r}\right) H_{r}\left( z_{r}\right) }f\left( z_{1},\dots,z_{n}\right)
dz_{1}\dots dz_{n}  \\
&=&\frac{1}{\left( 2\pi i\right) ^{n}}\int_{C_{1}}\dots\int_{C_{n}}\frac{1}{%
\prod_{r=1}^{n}\left( z_{r}-x_{r}\right) }f\left(
z_{1},\dots,z_{n}\right) dz_{1}\dots dz_{n}-\sum_{\mathbf{a\in A}}\sum_{\mathbf{%
k\in }\left[ \mathbf{0,}\nu \left( \mathbf{a}\right) -\mathbf{1}\right] }H_{%
\mathbf{a}}^{\mathbf{k}}\partial _{\mathbf{a}}^{\mathbf{k}}f \nonumber.
\end{eqnarray}

By applying the residue Theorem $n$-times, we obtain%
\begin{equation}\label{213}
\frac{1}{\left( 2\pi i\right) ^{n}}\int_{C_{1}}\dots\int_{C_{n}}\frac{1}{%
\prod_{r=1}^{n}\left( z_{r}-x_{r}\right) }f\left(
z_{1},\dots,z_{n}\right) dz_{1}\dots dz_{n}=f\left( x_{1},\dots,x_{n}\right) . 
\end{equation}%
From relations (\ref{212}) and (\ref{213}) we get (\ref{211}).
\end{proof}

\begin{rem}
{\rm
Applying Theorem \ref{errform} for $n=2$ we obtain the error formula for the bivariate Hermite polynomial, as given in \cite[Section 5]{GHITV}.
}
\end{rem}

\begin{rem}
{\rm
Lemma \ref{lemErr} provides an alternative way of generating the Hermite polynomial that interpolates any function $f$, provided that the integral can be calculated. %\textcolor{black}{A relevant example will be given in the next section.}
}
\end{rem}

\begin{rem}
{\rm
Lemma \ref{lemErr} and Theorem \ref{errform} hold also for complex functions of complex variables and for finite subsets $A_1,\dots ,A_n$ of $\mathbb{C}$.}
\end{rem}

\begin{rem}\label{rem8}
{\rm
Let us denote by $H\left( D_{z_{1}}\times \dots \times D_{z_{n}}\right) $ the
ring of all holomorphic functions on $D_{z_{1}}\times \dots\times D_{z_{n}}$
and by $\left( H_{1},\dots,H_{n}\right) $ the ideal generated by the
polynomials $H_{1},\dots,H_{n}$. Consider the linear operator $T_{\left( 
\mathbf{A,}\nu \right) }:H\left( D_{z_{1}}\times \dots \times D_{z_{n}}\right)
\rightarrow $ $V\left( \mathbf{A,}\nu \right) $ given by 
\begin{equation*}
T_{\left( \mathbf{A,}\nu \right) }\left( f\right) =\frac{1}{\left( 2\pi
i\right) ^{n}}\int_{C_{1}}\dots\int_{C_{n}}\frac{\prod_{r=1}^{n}\left(
H_{r}\left( z_{r}\right) -H_{r}\left( x_{r}\right) \right) }{%
\prod_{r=1}^{n}\left( z_{r}-x_{r}\right) H_{r}\left( z_{r}\right) }f\left(
z_{1},\dots,z_{n}\right) dz_{1}\dots dz_{n}.
\end{equation*}%
Clearly by Lemma
\ref{lemErr} we have that $T_{\left( \mathbf{A,}\nu \right) }$ is an
epimorphism and \begin{equation*}
\ker T_{\left( \mathbf{A,}\nu \right) }=\left\{ g\in
H\left( D_{z_{1}}\times \dots \times D_{z_{n}}\right):\,\partial _{\mathbf{a}}^{\mathbf{m}%
}g=0,~\mathbf{a\in A},~\mathbf{m\in }\left[ \mathbf{0,}\nu \left( \mathbf{a}%
\right) -\mathbf{1}\right] \right\} .
\end{equation*}
By using Leibniz derivative rule for
multivariable functions we have that for any $f\in H\left( D_{z_{1}}\times
\dots\times D_{z_{n}}\right) $ and any $g\in \ker T_{\left( \mathbf{A,}\nu
\right) }$ hold $\partial _{\mathbf{a}}^{\mathbf{m}}(fg)=0$. Therefore $\ker
T_{\left( \mathbf{A,}\nu \right) }$ is an ideal of $H\left( D_{z_{1}}\times
\dots\times D_{z_{n}}\right)$. In other words, $T_{\left( \mathbf{A,}\nu
\right) }$ is an ideal projector. Further by residue Theorem we get $
T_{\left( \mathbf{A,}\nu \right) }\left( H_{r}\right) =0,$ $r=1,\dots,n$ .
Therefore $\left( H_{1},\dots,H_{n}\right) \subseteq \ker T_{\left( \mathbf{A,}
\nu \right) }$.
}
\end{rem}

 In the next section the
connection between the ideals $\left( H_{1},\dots,H_{n}\right) $ and $\ker
T_{\left( \mathbf{A,}\nu \right) }$ will also be determined for any infinite field.

\section{On the ideal of the \textcolor{black}{proposed Hermite coordinate} interpolation}
\label{idealco}
\textcolor{black}{Finding the basis of the ideal of an interpolation is useful for the calculation of the interpolation error, as well as many other applications of algebraic geometry. This is a well-studied issue, with thorough theoretical investigation provided in \cite{deBoor2005} and \cite{ideal_Shekhtman} and more specific aspects in section 6 of \cite{gascasa}. In \cite{Sauer2018} Sauer constructs appropriate orthogonal $H$-basis of ideal interpolation which are involved in Prony's method to solve multidimensional Prony's problems.}

\textcolor{black}{Most of the work concentrates on Hermite interpolation of the total degree. In this section, we focus on the proposed coordinate Hermite interpolation.}
Let $A_{i},\;i=1,\dots,n,$ be finite subsets of a field $\Bbbk $ \textcolor{black}{ of characteristic zero}, and $\nu _{i}\left( a\right) $ the multiplicity of $a\in A_{i}$. Consider the polynomials $$
H_{i}\left( x_{i}\right) =\prod_{a\in A_{i}}\left(
x_{i}-a\right) ^{\nu _{i}\left( a\right) }\in \Bbbk \left[ x_{i}\right]
,\;i=1,\dots,n.$$

For $g\in \Bbbk \left[ x_{1},\dots,x_{n}\right] $, denote by\ $
\widetilde{g}$ the residue class of $g$ in the quotient ring $\Bbbk \left[
x_{1},\dots,x_{n}\right] /\left( H_{1},\dots,H_{n}\right) $. It is easy to verify that the linear operator $T_{\left( \mathbf{A,}\nu
,\Bbbk \right) }:\Bbbk \left[ x_{1},\dots,x_{n}\right] \rightarrow $ $V\left( 
\mathbf{A,}\nu \right) $, given by 
\begin{equation*}
T_{\left( \mathbf{A,}\nu ,\Bbbk \right) }\left( f\right) =\sum_{\mathbf{a\in
A}}\sum_{\mathbf{k\in }\left[ \mathbf{0,}\nu \left( \mathbf{a}\right) -%
\mathbf{1}\right] }H_{\mathbf{a}}^{\mathbf{k}}\partial _{\mathbf{a}}^{%
\mathbf{k}}f,
\end{equation*}%
is a projector. Further, similarly to Remark $8,$ it can be shown that $\ker
T_{\left( \mathbf{A,}\nu ,\Bbbk \right) }$ is an ideal of $\Bbbk \left[
x_{1},\dots,x_{n}\right] $. So $T_{\left( \mathbf{A,}\nu ,\Bbbk \right) }$ is
an ideal projector. \textcolor{black}{The following Lemma is well-known, however it is stated for completeness.}

\begin{lem}
The $\Bbbk$-algebra $\Bbbk \left[ x_{1},\dots,x_{n}\right] /\ker T_{\left( 
\mathbf{A,}\nu ,\Bbbk \right) }$ as a $\Bbbk$-vector space is isomorphic to 
$V\left( \mathbf{A,}\nu \right) $, and the following is valid 
\begin{equation*}
\dim _{\Bbbk }\Bbbk \left[ x_{1},\dots,x_{n}\right] /\ker T_{\left( \mathbf{A,}%
\nu ,\Bbbk \right) }=\prod_{i=1}^{n}\sum_{a_{i}\in A_{i}}\nu _{i}\left(
a\right)
\end{equation*}
\end{lem}
\begin{proof}
Let $\phi :\Bbbk \left[ x_{1},\dots,x_{n}\right] \rightarrow \Bbbk \left[
x_{1},\dots,x_{n}\right] /\ker T_{\left( \mathbf{A,}\nu ,\Bbbk \right) }$be
the canonical map. Consider any basis $B$ of the vector space $V\left( 
\mathbf{A,}\nu \right) .$ It is sufficient to show, that $\phi \left(
B\right) $ is a basis of $\Bbbk \left[ x_{1},\dots,x_{n}\right] /\ker
T_{\left( \mathbf{A,}\nu ,\Bbbk \right) }:$ Starting from $\sum_{b\mathbf{%
\in }B}c_{b}\phi \left( b\right) =0$, $c_{b}\in \Bbbk $ we have $\sum_{b%
\mathbf{\in }B}c_{b}\phi \left( b\right) \in \ker T_{\left( \mathbf{A,}\nu
,\Bbbk \right) }$, from which we obtain 
\begin{equation*}
\partial _{\mathbf{a}}^{\mathbf{m}}\sum_{b\mathbf{\in }B}c_{b}b=0,~\mathbf{%
a\in A}\text{ and}~\mathbf{m\in }\left[ \mathbf{0,}\nu \left( \mathbf{a}%
\right) -\mathbf{1}\right] .
\end{equation*}%
Further \textcolor{black}{ $~$because $B\subset V\left( \mathbf{A,}\nu \right) $ we have$\
\sum_{b\mathbf{\in }B}c_{b}b\in V\left( \mathbf{A,}\nu \right). $%
 Consequently by Corollary 2.4, the following holds:} $\sum_{b\mathbf{\in }B}c_{b}b=0.$
 Therefore, for each $b\mathbf{\in }B$, it holds that $c_{b}=0,$which means $\phi \left(
B\right) $ is linearly independent. Now we will show that $\phi \left(
B\right) $ generate $\Bbbk \left[ x_{1},\dots,x_{n}\right] /\ker T_{\left( 
\mathbf{A,}\nu ,\Bbbk \right) }$: By Theorem 2.2, for any $g\in \Bbbk \left[
x_{1},\dots,x_{n}\right] $ there is exactly one $f=\sum_{b\mathbf{\in }%
B}c_{b}b\in V\left( \mathbf{A,}\nu \right) $ such that $\partial _{\mathbf{a}%
}^{\mathbf{m}}f=\partial _{\mathbf{a}}^{\mathbf{m}}g,$ $\mathbf{a\in A},~%
\mathbf{m\in }\left[ \mathbf{0,}\nu \left( \mathbf{a}\right) -\mathbf{1}%
\right] .$ Consequently $\phi \left( g\right) =\phi \left( f\right) $, which
can be written as 
\begin{equation*}
\phi \left( g\right) =\sum_{b\mathbf{\in }B}c_{b}\phi \left( b\right) .
\end{equation*}
\end{proof}

\begin{lem}\label{lem42}
Given the polynomials $f_{i}\in \Bbbk \left[ x_{i}\right] ,~i=1,\dots,n~$, 
the quotient ring $\Bbbk \left[ x_{1},x_{n}\right] /\left(
f_{1},\dots,f_{n}\right) $ as a $\Bbbk$-vector space is isomorphic to $%
\left\{ f\in 
%TCIMACRO{\U{2102} }%
%BeginExpansion
\mathbb{C}
%EndExpansion
\left[ x_{1},\dots,x_{n}\right] :\deg _{i}f<\deg f_{i},i=1,\dots,n\right\} ,$
and 
\begin{equation*}
\dim _{\Bbbk }\Bbbk \left[ x_{1},\dots,x_{n}\right] /\left(
f_{1},\dots,f_{n}\right) =\prod_{i=1}^{n}\deg f_{i}.
\end{equation*}
\end{lem}

\begin{proof}
It is sufficient to prove that for any polynomiall $f\in \Bbbk %
\left[ x_{1},\dots,x_{n}\right] $ there exist unique $c_{\left(
k_{1},\dots ,k_{n}\right) }\in \Bbbk ,~0\leq k_{1}<\deg f_{1},\dots ,0\leq
k_{n}<\deg f_{n}$, such that 
\begin{equation*}
f(x_1,\dots, x_n)=\sum_{k_{1}=0}^{\deg f_{1}-1}\dots\sum_{k_{n}=0}^{\deg f_{n}-1}c_{\left(
k_{1},\dots ,k_{n}\right) }x_{1}^{k_{1}}\dots x_{n}^{k_{n}}+\left(
f_{1}(x_1),\dots,f_{n}(x_n)\right) .
\end{equation*}
To this end, we will apply the method of induction over $n$: For $n=1$
obviously the claim is valid. Let the claim be valid for $n-1,n>1$. Dividing $f$ $\ $by $f_{n}\left( x_{n}\right) $, we have $f=pf_{n}+q,$ where 
$p,q\in \Bbbk \left[ x_{1},\dots,x_{n}\right] $ are uniquely deteminated
holding $\deg _{n}q<\deg f_{n}$. Therefore
\begin{equation}\label{41}
f=\sum_{i=0}^{\deg f_{n}-1}q_{i}x_{n}^{i}+\left( f_{n}\right) , 
\end{equation}%
where $q_{i}$ are uniquely deteminated polynomials of $\Bbbk \left[
x_{1},\dots,x_{n-1}\right] $, and $\left( f_{n}\right) $ is the ideal of $%
\Bbbk \left[ x_{1},\dots,x_{n}\right] $ generated by $f_{n}$. Now by induction,
the polynomials $q_{i}$ can be written as, 
\begin{equation}\label{42}
q_{i}=\sum_{0\leq k_{1}<\deg f_{1}}\dots\sum_{0\leq k_{n-1}<\deg
f_{n-1}}c_{\left( k_{1},\dots,k_{n-1}\right)
}^{i}x_{1}^{k_{1}}\dots x_{n-1}^{k_{n-1}}+\left( f_{1},\dots,f_{n-1}\right) , 
\end{equation}%
where $c_{\left( k_{1},\dots,k_{n-1}\right) }^{i}\in \Bbbk $ are uniquely determined. By substitution of (\ref{42}) in (\ref{41}) we obtain the required expression.
\end{proof}

\begin{thm}\label{thm43}
The following holds $$\left( H_{1},\dots,H_{n}\right) =\ker T_{\left( \mathbf{A,}\nu
,\Bbbk \right) }.$$
\end{thm}

\begin{proof}
Applying the Leibniz derivative rule we conclude that for each $g\in \Bbbk \left[
x_{1},\dots,x_{n}\right] $ the following is valid $\partial _{\mathbf{a}}^{\mathbf{m}}\left(
gH_{\iota }\right) =0,~i=1,\dots,n\,\ $, for each $\mathbf{a\in A},~\mathbf{%
m\in }\left[ \mathbf{0,}\nu \left( \mathbf{a}\right) -\mathbf{1}\right] .$
Therefore $H_{i}\in \ker T_{\left( \mathbf{A,}\nu ,\Bbbk \right)
},~i=1,\dots,n$, or equivalently 
\begin{equation*}
\left( H_{1},\dots,H_{n}\right) \subseteq \ker T_{\left( \mathbf{A,}\nu ,\Bbbk
\right) }.
\end{equation*}%
Now it is sufficient to show that $\ker T_{\left( \mathbf{A,}\nu ,\Bbbk \right)
}\subseteq \left( H_{1},\dots,H_{n}\right) $: Let $g$ be any element of $I$. By
Lemma \ref{lem42} there is unique polynomial $p\in \left\{ f\in \Bbbk \left[
x_{1},\dots,x_{n}\right] :\deg _{i}f<\deg H_{i},i=1,\dots,n\right\} =V\left( 
\mathbf{A,}\nu \right) $ such that 
\begin{equation}\label{43}
g=p+\left( H_{1},\dots,H_{n}\right)  
\end{equation}%
From (\ref{43}) we obtain that for any $\mathbf{m\in }\left[ \mathbf{%
0,}\nu \left( \mathbf{a}\right) -\mathbf{1}\right] $ and any$~\mathbf{a\in A}
$ 
\begin{equation*}
\partial _{\mathbf{a}}^{\mathbf{m}}g=\partial _{\mathbf{a}}^{\mathbf{m}}p,
\end{equation*}%
which implies $\partial _{\mathbf{a}}^{\mathbf{m}}p=0$. Therefore, 
Corollary 2.4 \ yields that $p=0,$ which combined with $\left( 4.3\right) $
gives $g\in \left( H_{1},\dots,H_{n}\right) .$
\end{proof}

\begin{rem}
By Theorem \ref{thm43} it follows that $\left\{ H_{1}\left( x_{1}\right)
,\dots,H_{n}\left( x_{n}\right) \right\} $ is a reduced Gr{\"o}bner basis of $\ker T_{\left( \mathbf{A,}\nu ,\Bbbk \right) }$.
\end{rem}

Given the polynomials $f_{i}\in \Bbbk \left[ x_{i}\right] ,~i=1,\dots,n$, $A_{i}$, the sets of their roots in the algebraic closure $\overline{\Bbbk }$
of $\Bbbk $, and $\nu _{i}\left( a\right) $ the multiplicity of $a\in A_{i}$, we obtain $f_{i}\left( x_{i}\right) =\prod_{a\in A_{i}}\left( x_{i}-a\right)
^{\nu _{i}\left( a\right) },~i=1,\dots,n$. The ideal of $\overline{\Bbbk }\left[
x_{1},\dots,x_{n}\right] $ is denoted by $\overline{I}$ and it is   generated by $f_{1},\dots,f_{n}$. For $g\in \Bbbk \left[ x_{1},\dots,x_{n}\right] $ let us denote by $\overline{g%
}\ $the residue class of $g$\ in $\overline{\Bbbk }\left[ x_{1},\dots,x_{n}%
\right] /\overline{I}$ .

\begin{cor}\label{cor44}
The ideal $\left( f_{1},\dots,f_{n}\right) $ of $\Bbbk \left[ x_{1},\dots,x_{n}%
\right] $ is equal to $\ker T_{\left( \mathbf{A,}\nu ,\overline{\Bbbk }%
\right) }\cap \Bbbk \left[ x_{1},\dots,x_{n}\right] .$
\end{cor}

\begin{proof}
By Theorem \ref{thm43} we get $\left( f_{1},\dots,f_{n}\right) \subseteq
\ker T_{\left( \mathbf{A,}\nu ,\overline{\Bbbk }\right) }\cap \Bbbk \left[
x_{1},\dots,x_{n}\right] $. Therefore, it is sufficient to prove that $\ker
T_{\left( \mathbf{A,}\nu ,\overline{\Bbbk }\right) }\cap \Bbbk \left[
x_{1},\dots,x_{n}\right] \subseteq \left( f_{1},\dots,f_{n}\right) :$ \ Let $f$
be any element of $J\left( \mathbf{A,}\nu \right) \cap \Bbbk \left[
x_{1},\dots,x_{n}\right] .$ Then by Theorem \ref{thm43} we get
\begin{equation}\label{44}
f=\sum_{i=1}^{n}f_{i}g_{i}, 
\end{equation}%
where $g_{i}\in \overline{\Bbbk }\left[ x_{1},\dots,x_{n}\right]$. Consider
the finite field extension $\Bbbk \subset \Bbbk_{1}$, where $\Bbbk_{1}$ is
the field constructed by adjunction of all coefficients of the polynomials $g_{i}.$ Let $r$ be the degree of $\Bbbk \subset \Bbbk _{1}$ and $\rho $ a
primitive element. Then (\ref{44}) can be written as 
\begin{equation}\label{45}
f=\sum_{j=0}^{r-1}\rho ^{j}\sum_{i=1}^{n}f_{i}g_{ij},  
\end{equation}
where $g_{ij}\in \Bbbk \left[ x_{1},\dots,x_{n}\right]$. Finally, since $
\left\{ 1,\rho ,\dots,\rho ^{r-1}\right\} $ is $\Bbbk$- linearly independent, from (\ref{45}) we obtain
\begin{equation*}
f=\sum_{i=1}^{n}f_{i}g_{i0}.
\end{equation*}%
That is $f\in I.$
\end{proof}

\begin{rem}\label{rem10}
Corollary \ref{cor44} can also be formulated  as: $f\in \left(
f_{1},\dots,f_{n}\right) \subset \Bbbk \left[ x_{1},\dots,x_{n}\right] $, if and
only if $ \partial _{\mathbf{a}}^{\mathbf{m}}f=0,$ for each $\mathbf{a\in A}
$ and$~\mathbf{m\in }\left[ \mathbf{0,}\nu \left( \mathbf{a}\right) -\mathbf{%
1}\right]$.
\end{rem}

%%%%%%%%%%%%%%%%%%%%%%%%%%%%%%%%%%%%%%%%
%%%%%%%%%%%%%%%%%%%%%%%%%%%%%%%%%%%%%%%
Let $A_{i},~i=1,\dots,n,$ be finite subsets of an infinite field $\Bbbk $, and $
\nu _{i}\left( a\right) $ the multiplicity of $a\in A_{i}$, and  the
polynomials $H_{i}\left( x_{i}\right) =\prod_{a\in A_{i}}\left(
x_{i}-a\right) ^{\nu _{i}\left( a\right) }\in \Bbbk \left[ x_{i}\right]
,~i=1,\dots,n.$ Given any polynomial $g\in \Bbbk \left[ x_{1},\dots,x_{n}\right] 
$ we apply the Euclidean polynomial division $n$-times.
Assuming that we start with the 1st dimension, there are unique polynomials $q_{1},r_{1}\in \Bbbk \left[
x_{1},\dots,x_{n}\right] $ with $\deg _{1}r_{1}<\deg H_{1}$, such that $%
g=H_{1}q_{1}+r_{1},$. Continuing with the 2nd dimension ($i=2$), there are unique polynomials $q_{2},r_{2}\in \Bbbk %
\left[ x_{1},\dots,x_{n}\right] $ with $\deg _{i}r_{2}<\deg H_{i}$,~$i=1,2~$%
such that $r_{1}=H_{2}q_{2}+r_{2}.$ Continuing in a similar manner until the last dimension, $n$, we obtain unique polynomials $%
q_{n},r_{n}\in \Bbbk \left[ x_{1},\dots,x_{n}\right] $ with $\deg
_{i}r_{n}<\deg H_{i}$,~$i=1,2,\dots,n~$such that $r_{n-1}=H_{n}q_{n}+r_{n}.$
Combining the above $\ n$ polynomial divisions we get 
\begin{equation}\label{rel46}
g=r_{n}+\sum_{i=1}^{n}H_{i}q_{i}.  
\end{equation}

It is easy to verify that for $n$ dimensions there are $n!$ different orders of application of the polynomial divisions, all of which lead to the same result $r_n.$  

\begin{prop}\label{proer1}
Given any polynomial $g\in \Bbbk \left[ x_{1},\dots,x_{n}\right] $, the
interpolating polynomial of the interpolation as presented in Theorem
\ref{thm1} with $t_{\mathbf{a}}^{\mathbf{m}}:=\partial _{\mathbf{a}}^{\mathbf{m}}g,~%
\mathbf{a\in A},~\mathbf{m\in }\left[ \mathbf{0,}\nu \left( \mathbf{a}%
\right) -\mathbf{1}\right] $ is equal to $r_{n}$, and the remainder is equal
to $$\sum_{i=1}^{n}H_{i}q_{i}.$$

\end{prop}

\begin{proof}
From $\deg _{i}r_{n}<\deg H_{i}$,~$i=1,2,\dots,n~$\ we have that $r_{n}\in
V\left( \mathbf{A,}\nu \right) ,$ and from (\ref{rel46}) by Remark \ref{rem10} we obtain $\partial _{\mathbf{a}}^{\mathbf{m}}g=\partial _{\mathbf{a}}^{%
\mathbf{m}}r_{n},~\mathbf{a\in A},~\mathbf{m\in }\left[ \mathbf{0,}\nu
\left( \mathbf{a}\right) -\mathbf{1}\right]$. These arguments in combination
with (\ref{rel46}) lead us immediately to the claims of the current
proposition.
\end{proof}

Proposition \ref{proer1} can be used to express the remainder of the Hermite interpolation in closed form using the cascaded polynomial divisions in case of interpolating a polynomial function. In the general case of any differentiable function to be interpolated (other than polynomial), this error calculation is applicable after expanding the function using Taylor Theorem. The following examples demonstrates error calculation for a univariate and a 3-D polynomial function.

\begin{exam} \label{hermerrr}
{\rm
Let $g=(x-2)^{11}+x^4+9$. Consider the grid $A=\{0.7, 1.2 ,1.7, 2.2\}$, with multiplicity function $\nu(a) = 2 $, \textcolor{black}{for every point $a$ on the grid $ A$}. We derive the polynomial $H_1(x)= (x-0.7)^2(x-1.2)^2(x-1.7)^2(x-2.2)^2$.
Applying Eq.(\ref{rel46}), we obtain the quotient \textcolor{black}{$q=x^3-\frac{52}{5}x^2+\frac{2087}{50}x-\frac{1918}{25}$} and \textcolor{black}{the remainder $r$ of the polynomial division} $\frac{g}{H_1}.$
According to Eq.(\ref{rel46}), $r$ is the Hermite polynomial with \textcolor{black}{$ \deg r\leq \sum_i \left(\nu \left(a_i\right)\right)-1=7$}, given by

\begin{eqnarray*}
r = \frac{9273}{178} x^7   
    -593 x^6    +      \frac{31499}{11} x^5     -\frac{106481}{14}x^4 +      \frac{108109}{9} x^3      -\frac{349850}{31}x^2    
   +\frac{40932}{7}x -\frac{6409}{5} .
\end{eqnarray*}
The error of the Hermite interpolation is easily obtained by
\begin{eqnarray*}
R &=&H_1 q_1=
x^{11} -22x^{10}   +220x^9  -1320x^8     +\frac{52279}{10}x^7
  -14191x^6 +\frac{293749}{11}x^5 -\frac{173166}{5}x^4 +\frac{272051}{9}x^3\\ &&-\frac{33749}{2}x^2 
   +\frac{37916}{7}x -\frac{3786}{5}.
\end{eqnarray*}

It can be observed that the Hermite polynomial $r(x)$ is unique among the polynomials of degree up to 7 that satisfies the requirements $g(a_i)=r(a_i) \text{ and }\frac{dg}{dx}=\frac{d r}{d x}$ for each supporting point $a_i$. The interpolated function $g(x)$ and the resulting Hermite polynomial are depicted in Figure \ref{hermer} (a). The interpolation error is plotted in the semi-logarithmic graph of Figure \ref{hermer} (b). The error remains very low in $[a_1,a_4]$ and approaches 0 very fast in the vicinity of the $a_i$.

%%%%%%%%%%%%%%%%%%%%%
\begin{exam}
Let us consider $g=(x-2y)^7 +(x+z^2)^4$ and the grid $\mathbf{A}=\{-1, 0 ,1\}\times\{-1, 0 ,1\}\times\{-1, 0 ,1\}$, with multiplicity function $\nu(\mathbf{a}) = (2,2,2) $, for every point on the grid $\mathbf{a} \in \mathbf{A}$. We derive the polynomials \textcolor{black}{$H_1(x)= (x+1)^2x^2(x-1)^2$}, \textcolor{black}{$H_2(y)= (y+1)^2 y^2(y-1)^2$} and \textcolor{black}{$H_3(z)= (z+1)^2 z^2(z-1)^2$}.
Applying Eq.(\ref{rel46}), we obtain the quotient \textcolor{black}{$q_1=x-14y$} and the remainder $r_1$ of the polynomial division $\frac{g}{H_1}.$, Subsequently, we obtain the quotient \textcolor{black}{$q_2=448x-128y$} and the remainder $r_2$ of the polynomial division $\frac{r_1}{H_2}$ and finally, the quotient \textcolor{black}{$q_3=z^2+4x+2$} and the remainder $r_3$ of the polynomial division $\frac{r2}{H_3}.$
According to Eq.(\ref{rel46}), $r_3$ is the Hermite polynomial given by
$$r_3=84x^5y^2 + 2x^5 - 280x^4 y^3 - 28x^4y + 560x^3y^4 - 672x^2y^5 + 14x^2y + 3x^2z^2 + 896x y^4 - 448 x y^2 + 3 x z^4 - 256 y^5 + 128 y^3 + 2 z^4 - z^2$$
The remainder $R=H_1q_1+H_2q_2+H_3q_3$ is calculated as:
$$R=z^2(z^2 - 1)^2 + x^2(x - 14y)(x - 1)^2(x + 1)^2 + y^2(448x - 128y)(y - 1)^2 (y + 1)^2.$$
\end{exam}

%%%%%%%%%%%%%%%%%%%%%%

\begin{figure}[h!]
\centering
\subfigure[]{{\includegraphics[width=0.6\textwidth]{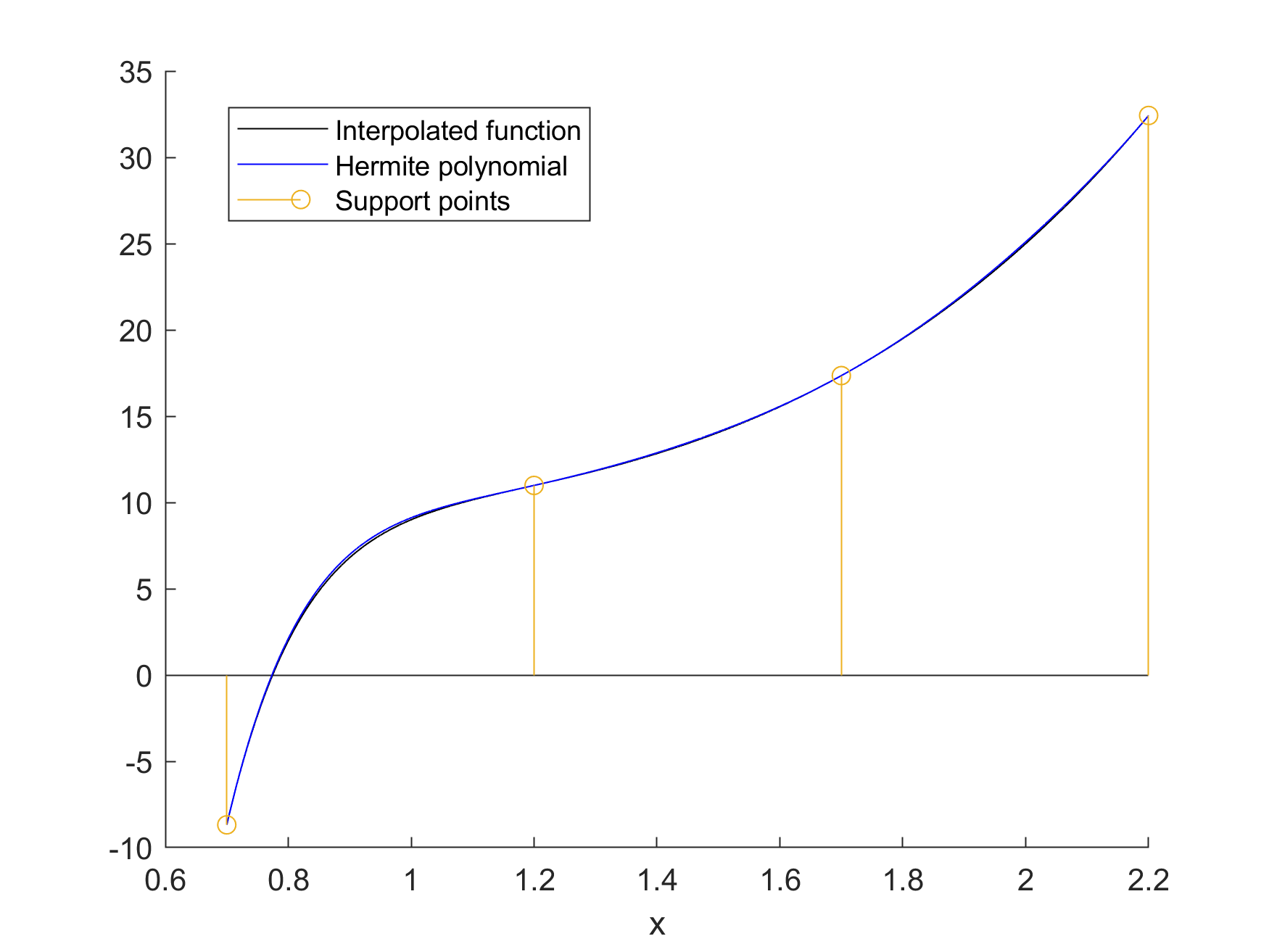}}}\hfill
\subfigure[]{{\includegraphics[width=0.6\textwidth]{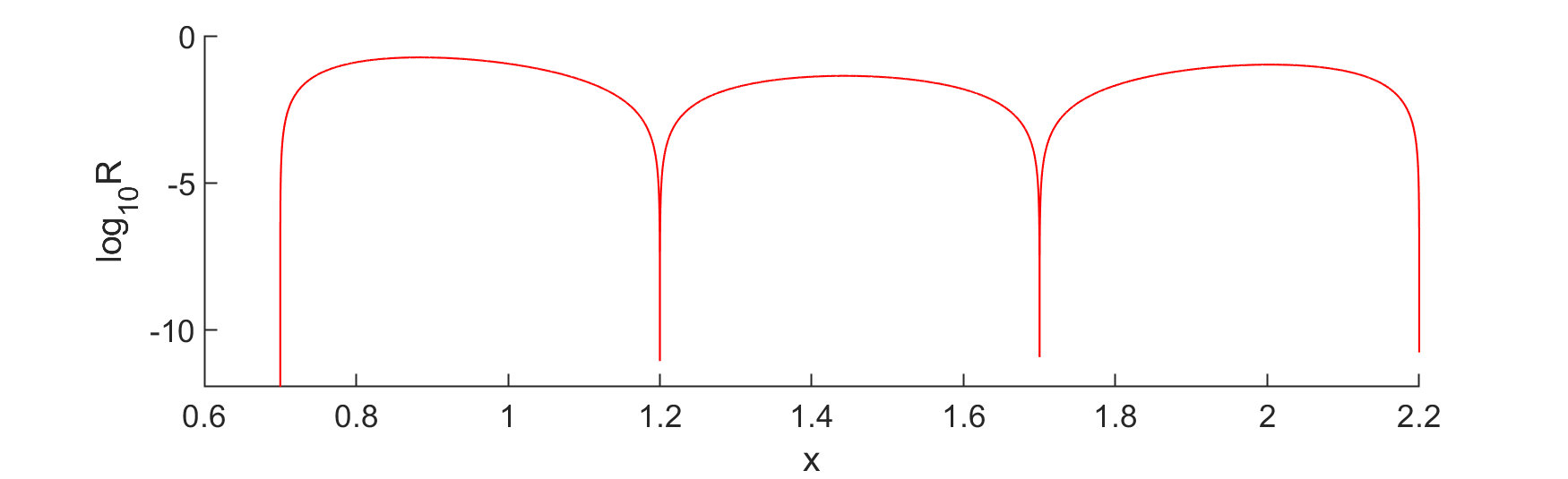}}}
\caption{
(a) The interpolated function and Hermite polynomial of the Example \ref{hermerrr}, along with the support points  (b) Semi-logarithmic plot of the remainder of interpolation.}
\label{hermer}
\end{figure}

}
\end{exam}

%Previous Section
%\input{idealconn}

%\newpage
\section{An alternative quick proof of Theorem \ref{th1}}
\label{altproof}
In this section we provide a quick proof of Theorem \ref{th1}, as an alternative to the proof given in \cite{fractmult}. 
In addition to the already defined $H_{a}^{k}(x)$, we define the following

%    \item Let $A$ be a finite subset of $\Bbbk $. Then for any $a\in A$ and $%    k\in \left\{ 1,\dots,\nu \left( a\right) -1\right\} $,     \begin{equation*}    H_{a}^{k}\left( x\right) :=H_{a}\left( x\right) \frac{\left( x-a\right) ^{k}%    }{k!}\sum_{t=0}^{\nu \left( a\right) -k-1}\frac{\left( x-a\right) ^{t}}{t!}    \left( \frac{1}{H_{a}\left( x\right) }\right) ^{\left( t\right) }\left(    a\right) \in \Bbbk \left[ x\right] ,    \end{equation*} where the superscript $(t)$ indicates order of derivative.

\begin{equation*}
H_{\mathbf{a}}^{\mathbf{k}}\left( x_{1},\dots, x_{n}\right)
:=\prod_{i=1}^{n}H_{a_{i}}^{k_{i}}\left( x_{i}\right) \in \Bbbk \left[
x_{1},\dots,x_{n}\right], \text{ and}
\end{equation*}

%Let us define $Q\left( \mathbf{A,}\nu \right)$ as follows

\begin{equation*}
Q\left( \mathbf{A,}\nu \right) :=\left\{ H_{\mathbf{a}}^{\mathbf{k}}:\mathbf{%
a\in A}\text{, }\mathbf{k\in }\left[ \mathbf{0,}\nu \left( \mathbf{a}\right)
-\mathbf{1}\right] \right\} .
\end{equation*}

First, we state the following remarks and lemma.
%%%%%%%%%%%%%%%%%%%%%%%%%%%%%5
\begin{rem}\label{rem9}
{\rm
From $\deg _{i}H_{\mathbf{a}}^{\mathbf{k}}=k_{i}+\sum_{a\in A_{i}-\left\{
a_{i}\right\} }\nu _{i}\left( a\right) <\sum_{a\in A_{i}}\nu _{i}\left(
a\right) $ follows $Q\left( \mathbf{A,}\nu \right) \subset V\left( \mathbf{A,%
}\nu \right) $, and $\left\vert Q\left( \mathbf{A,}\nu \right) \right\vert
=\prod_{i=1}^{n}\sum_{a_{i}\in A_{i}}\nu _{i}\left( a\right) =\dim
_{k}V\left( \mathbf{A,}\nu \right) .$}
\end{rem}

\begin{rem}\label{remkron}
{\rm
For any $a,~b\in A$, and $k,m\in \left\{ 1,\dots,\nu \left( a\right)
-1\right\} $ we have (see \cite{SPITZ})
\begin{equation*}
\partial ^{m}H_{a}^{k}\left( b\right) =\left\{ 
\begin{array}{ccc}
0, & \text{if} & a\neq b \\ 
\delta _{km} & \text{if} & a=b%
\end{array}%
\right. ,
\end{equation*}%
where $\delta _{km}$ is the Kronecker Delta. From this we immediately
obtain, for any $\mathbf{a,b\in A}$, and any $~\mathbf{k,~m\in }\left[ 
\mathbf{0,}\nu \left( \mathbf{a}\right) -\mathbf{1}\right] $ we have%
\begin{equation*}
\partial _{\mathbf{a}}^{\mathbf{m}}H^{\mathbf{k}}\left( \mathbf{b}\right)
=\left\{ 
\begin{array}{ccc}
0, & \text{if} & \mathbf{a}\neq \mathbf{b} \\ 
\delta _{\mathbf{km}} & \text{if} & \mathbf{a}=\mathbf{b}%
\end{array}%
\right. .
\end{equation*}}
\end{rem}

\begin{lem}\label{lemex}
$Q\left( \mathbf{A,}\nu \right) $ is a basis of $V\left( \mathbf{A,}\nu
\right)$.
\end{lem}

\begin{proof}
By Remark \ref{rem9} it is sufficient to show that $Q\left( \mathbf{A,}\nu
\right) $ is linearly independent: From $\sum_{\mathbf{b\in A}}\sum_{\mathbf{%
k\in }\left[ \mathbf{0,}\nu \left( \mathbf{b}\right) -\mathbf{1}\right] }c_{%
\mathbf{b}}^{\mathbf{m}}H_{\mathbf{b}}^{\mathbf{m}}=0,~c_{\mathbf{b}}^{%
\mathbf{m}}\in 
%TCIMACRO{\U{2102} }%
%BeginExpansion
\mathbb{C}
%EndExpansion
$ for each $\mathbf{a\in A,~k\in }\left[ \mathbf{0,}\nu \left( \mathbf{a}%
\right) -\mathbf{1}\right] $ follows $\sum_{\mathbf{b\in A}}\sum_{\mathbf{%
k\in }\left[ \mathbf{0,}\nu \left( \mathbf{b}\right) -\mathbf{1}\right] }c_{%
\mathbf{b}}^{\mathbf{m}}\partial _{\mathbf{a}}^{\mathbf{k}}H_{\mathbf{b}}^{%
\mathbf{m}}=0,$ from which by Remark \ref{remkron} we conclude $c_{\mathbf{a}}^{%
\mathbf{k}}=0$.
\end{proof}

\begin{thm}
Given $t_{\mathbf{a}}^{\mathbf{k}}\in \,,~\mathbf{a\in A},~\mathbf{k\in }%
\left[ \mathbf{0,}\nu \left( \mathbf{a}\right) -\mathbf{1}\right] .$ Then
there is a unique $p\in V\left( \mathbf{A,}\nu \right) $ such that $\partial
_{\mathbf{a}}^{\mathbf{k}}p=t_{\mathbf{a}}^{\mathbf{k}},~\mathbf{a\in A}, $ $
\mathbf{k\in }\left[ \mathbf{0,}\nu \left( \mathbf{a}\right) -\mathbf{1}%
\right] $, and $p$ is given by 
\begin{equation*}
p=\sum_{\mathbf{a\in A}}\sum_{\mathbf{m\in }\left[ \mathbf{0,}\nu \left( 
\mathbf{a}\right) -\mathbf{1}\right] }t_{\mathbf{a}}^{\mathbf{k}}H_{\mathbf{a%
}}^{\mathbf{k}}=\sum_{a_{1}\in A_{1}}\dots\sum_{a_{n}\in
A_{n}}\sum_{k_{n}=0}^{\nu _{1}\left( a_{1}\right) -1}\dots\sum_{k_{n}=0}^{\nu
_{n}\left( a_{n}\right) -1}\ \prod_{j=1}^{n}\partial
_{i}^{k}H_{a_{j}}^{k_{j}}\left( x_{j}\right) t_{\mathbf{a}}^{\mathbf{k}}.
\end{equation*}
\end{thm}

\begin{proof}
By Remark \ref{remkron} for each $\mathbf{a\in A,~k\in }\left[ \mathbf{0,}\nu \left( 
\mathbf{a}\right) -\mathbf{1}\right] $ we have 
\begin{equation*}
\partial _{\mathbf{a}}^{\mathbf{k}}p=\sum_{\mathbf{b\in A}}\sum_{\mathbf{%
k\in }\left[ \mathbf{0,}\nu \left( \mathbf{b}\right) -\mathbf{1}\right] }t_{%
\mathbf{b}}^{\mathbf{m}}\partial _{\mathbf{a}}^{\mathbf{k}}H_{\mathbf{b}}^{%
\mathbf{m}}=t_{\mathbf{a}}^{\mathbf{k}}.
\end{equation*}
The uniqueness of $p$ follows\ By Lemma \ref{lemex}.
\end{proof}
%\newpage

%\section{Computational Complexity}
%\label{comp}
%\input{rema1}

\section{Hermite splines interpolation approach}
\label{hermitesplines}
In this section, we present the $n$-dimensional Hermit interpolation via splines, which results in a less computationally intensive implementation. As it can be observed from Eq. (\ref{Hmikra}) and (\ref{Hmegala}), the maximum degree of each of the univariate Hermite polynomials for each dimension is given by the following formula:
\begin{equation}
    %\max_{a\in A_i}\{v_i(a)\} -\min_{a\in A_i}\{v_i(a)\}- 1 + \sum_{a \in A_i} v_i(a) =
     \max_{a\in A_i}\min_{b \in A_i}\{v_i(a)-v_i(b)\} - 1 +  \sum_{a \in A_i} v_i(a).
\end{equation}

The maximum total degree of the polynomial is:
\begin{equation}
   \sum_{i=1}^{n}\left(\max_{a\in A_i}\min_{b \in A_i}\{v_i(a)-v_i(b)\} - 1 +  \sum_{a \in A_i} v_i(a) \right)
\end{equation}

\begin{figure}[!h]
     \centering
         \includegraphics[width=9cm]{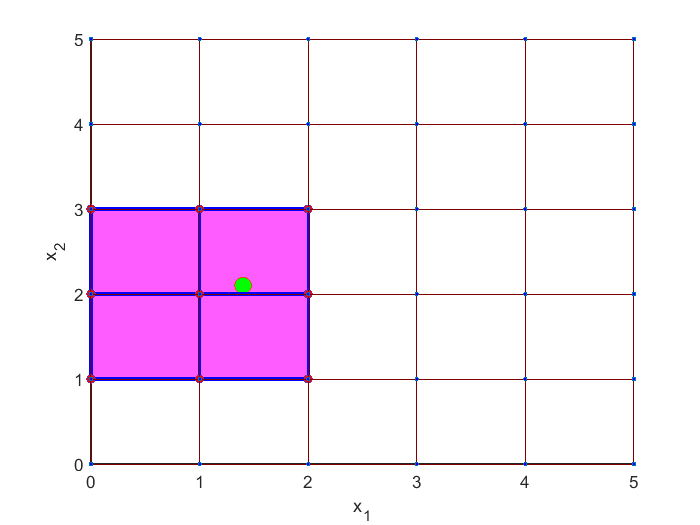}
         \caption{Visualization of an example of spline local support points round a random point (green circle) in 2D, with $\lvert B_1 \rvert=\lvert B_2 \rvert=3$. }
         \label{figsplinevis2d}
\end{figure}

\begin{figure}[!h]
     \centering
         \includegraphics[width=12cm]{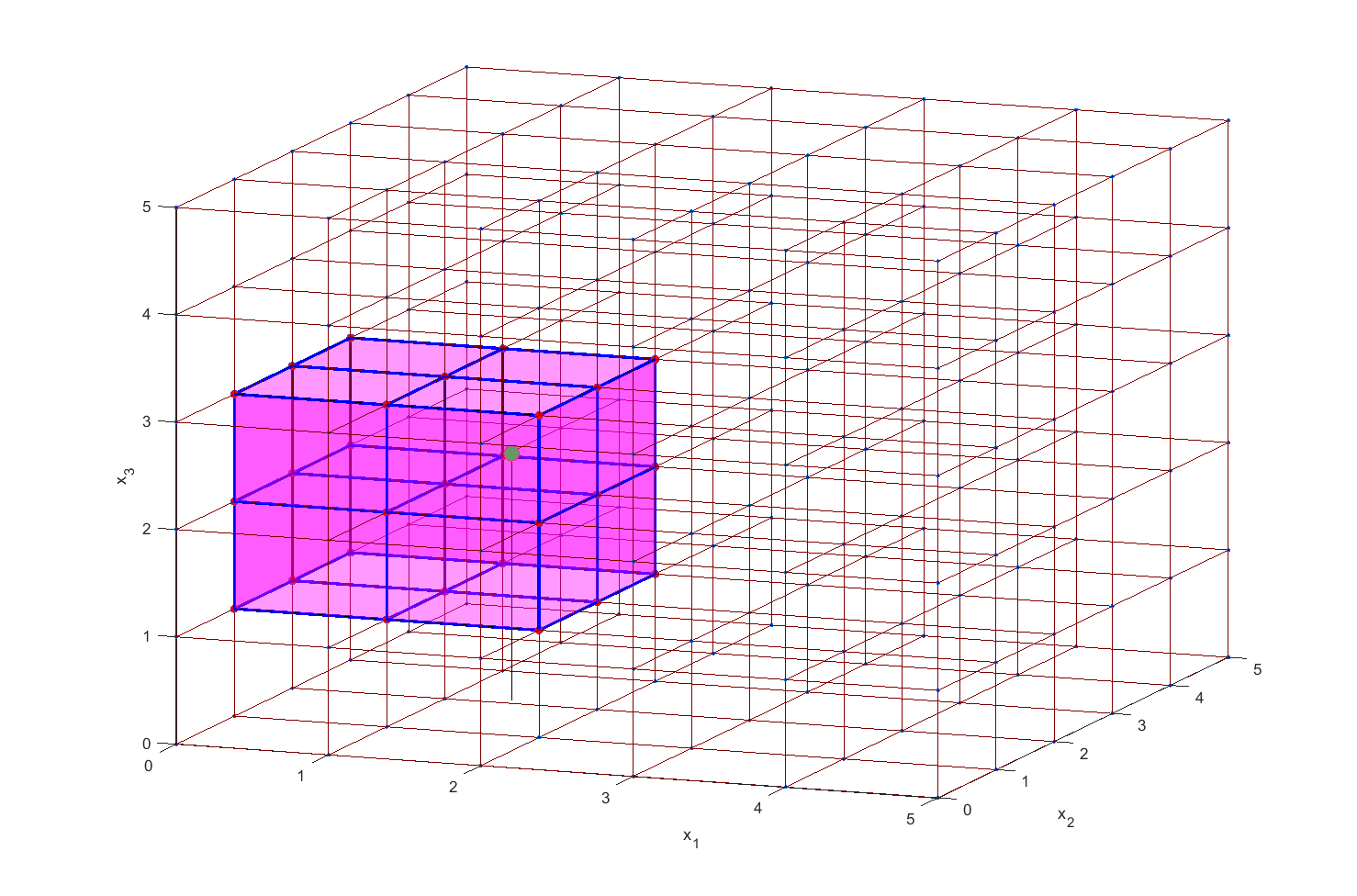}
         \caption{Visualization of an example of spline local support points round a random point (green circle) in 3D, with $\lvert B_1 \rvert=\lvert B_2 \rvert=\lvert B_3 \rvert=3$. }
         \label{figsplinevis3d}
\end{figure}

Thus, even for small numbers of support points on each axis and small values of multiplicity function $v_{i}$ the maximum total degree is of considerable value. For example if $|A_1|\times |A_2| \times |A_3| = 7\times 7 \times 7$, and $\mathbf{v} = (v_1(a),v_2(a),v_3(a))=(3,3,3)$, for every point on the grid $a \in \mathbf{A}$, the resulting  degree of the polynomials for each dimension is $7\cdot 3 +0-1=20$, and the total degree is $60$. Such high-degree polynomials are numerically unstable, since very small variations of the high-order coefficients may cause significant errors in the value of the polynomial. For this reason, the piece-wise calculation of the Hermite polynomial using only the local support points for each position $\bf x$ is preferred. More specifically, if the number of local support points is $|B_{i}|<|A_{i}|$ for each dimension $i$, then for any point $\mathbf{q}\in\left[\min\{A_{1}\},\max\{A_{1}\}\right]\times \left[\min\{A_{2}\},\max\{A_{2}\}\right]\times...\times\left[\min\{A_{n}\},\max\{A_{n} \}\right]$, 
the local support points along the $i^{\text{th}}$ dimension can be defined as follows. In the case of even $|B|$ the local support points along the $i$ dimension a defined as the closest $\frac{|B|}{2}$ points less than $q_i$ and the $\frac{|B|}{2}$ points closest and greater than $q_i$. 
In the case of odd $|B|$, the local support points along the $i$ dimension are defined as the single support point with $ith$ coordinate closest to $q_i$, and the $\frac{|B|- 1}{2}$ support points less than the closest and the $\frac{|B|- 1}{2}$ support points with $ith$ coordinate greater than the closest support point. 
In the most common case, where the Hermite interpolation is applied to signals in one or more dimensions, the support points are assumed to have integer coordinates, and the local support points are defined as follows:

\begin{equation}
   B_i:= \begin{cases}
     \left[ [q_{i}] -\frac{|B_i|-1}{2}, \dots, [q_{i}] +\frac{|B_i|-1}{2}\right]&,  \; \text{if }n\; |B_i|\mod 2=1 \\
     &\\
     \left[ \lfloor q_{i} \rfloor {-\frac{|B_i|}{2}+1}, \dots, \lfloor q_{i} \rfloor +\frac{|B_i|}{2}\right]&,  \; \text{otherwise}
   \end{cases}, \;i=1,\dots, n,
   \label{spline_local_sup}
\end{equation}
where $\lfloor .\rfloor$ denotes the rounding to the smallest nearest integer (floor) operator and $[.]$ denotes the rounding to the nearest integer (round) operator.
Fig. \ref{figsplinevis2d} and \ref{figsplinevis3d} visualize the support points for a random point $\mathbf{q}$ denoted by the green circle and $\lvert B_1 \rvert=\lvert B_2 \rvert=3$ in 2D and $\lvert B_1 \rvert=\lvert B_3 \rvert=\lvert B_3 \rvert=3$ in 3D, respectively, according to Eq.(\ref{figsplinevis2d}).
%3,4 masks 

%%%%%%%%%%%%%%%%%%%%%%%%%%%%%%%%%%%%%%%%%%

\subsection{Continuity of the proposed Hermite Splines}
\label{contsplines}
\textcolor{black}{
In this subsection, we provide formal proof that the proposed Hermite polynomials form a continuous spline. To this end, we will consider two $n$-dimensional grids, $\mathbf{A}_{a-}^{k}$ and $\mathbf{A}_{a+}^{k}$, which share a common hyperplane of support points normal to the $k^{th}$ dimension at position $a_k=a$. It will suffice to prove that the two correspondingly defined polynomials are continuous along the common hyperplane.
}

\textcolor{black}{
More formally, let $A_{i},~i=1,\dots,n$, be finite subsets of $
\mathbb{R}$, and $\nu _{i}\left( a\right) $ be the multiplicity of $a\in {A_i}$, and finally $\ t_{\mathbf{a}}^{\mathbf{k}}\in \Bbbk \,,~\mathbf{a\in A},$ $\mathbf{k\in }\left[ 
\mathbf{0,}\nu \left( \mathbf{a}\right) -\mathbf{1}\right]$. Consider the
polynomials $H_{i}\left( x_{i}\right) =\prod_{a\in A_{i}}\left(
x_{i}-a\right) ^{\nu _{i}\left( a\right) }\in \Bbbk \left[ x_{i}\right]
,~i=1,\dots,n$. Further for any $k=1,\dots,n$, any $a\in A_{k},$ 
let us denote $\mathbf{A}_{a-}^{k}:=\left\{ \mathbf{a\in A:}\text{ }a_{k}\leq a\right\} $, $%
\mathbf{A}_{a+}^{k}:=\left\{ \mathbf{a\in A:}\text{ }a_{k}\geq a\right\} $,
and $\mathbf{A}_{a}^{k}:=\left\{ \mathbf{a\in A:}\text{ }a_{k}=a\right\} .$
Clearly the following is valid $\mathbf{A}_{a-}^{k}\cap \mathbf{A}_{a+}^{k}=\mathbf{A}%
_{a}^{k}$. Consider the vector spaces $V\left( \mathbf{A}_{a-}^{k}\mathbf{,}%
\nu \right) $, $V\left( \mathbf{A}_{a+}^{k}\mathbf{,}\nu \right) \,$as
defined in the introduction and the vector space $V\left( \mathbf{A}_{a}^{k}%
\mathbf{,}\nu \right) $ consisting of all polynomials $f\left(
x_{1},\dots,x_{k-1},a,x_{k+1}\dots,.,x_{n}\right) $, such that $\deg
_{i}f<\sum_{a\in A_{i}}\nu _{i}\,\left( a\right) ,i=1,\dots,k-1,k+1,\dots,n.$ We
denote by $p_{a-}^{k}$ and $\ p_{a+}^{k}$  the following two unique Hermite interpolant polynomials (see Theorem
\ref{thm1}) of the vector spaces $V\left( \mathbf{A}_{a-}^{k}%
\mathbf{,}\nu \right) $ and $V\left( \mathbf{A}_{a+}^{k}\mathbf{,}\nu
\right) $ respectively, which satisfy the following conditions:%
\begin{equation}\label{eqtaged1}
\partial _{\mathbf{a}}^{\mathbf{m}}p_{a-}^{k}=t_{\mathbf{a}}^{\mathbf{m}},~%
\mathbf{a\in A}_{a-}^{k},~\mathbf{m\in }\left[ \mathbf{0,}\nu \left( \mathbf{%
a}\right) -\mathbf{1}\right],   
\end{equation}%
and 
\begin{equation}\label{eqtaged2}
\partial _{\mathbf{a}}^{\mathbf{m}}p_{a+}^{k}=t_{\mathbf{a}}^{\mathbf{m}},~%
\mathbf{a\in A}_{a+}^{k},~\mathbf{m\in }\left[ \mathbf{0,}\nu \left( \mathbf{%
a}\right) -\mathbf{1}\right] . 
\end{equation}
} 

\begin{thm}
\textcolor{black}{
    The Hermite polynomials $\left\{ 
\begin{array}{cc}
p_{a-}^{k}\left( \mathbf{x}\right) , & \mathbf{x\in A}_{a-}^{i} \\ 
p_{a+}^{k}\left( \mathbf{x}\right) , & \mathbf{x\in A}_{a+}^{i}%
\end{array}%
\right. $ form a continuous spline for any $\mathbf{x}$ with $x_k=a$.
}
\end{thm}

\begin{proof}
\textcolor{black}{
It is sufficient to prove the following identity: 
$$
p_{a-}^{k}\left( x_{1},\dots,x_{k-1},a,x_{k-1},\dots,x_{n}\right)
=p_{a+}^{k}\left( x_{1},\dots,x_{k-1},a,x_{k-1},\dots,x_{n}\right). 
$$
Given $p_{a-}^{k}\in V\left( \mathbf{A}_{a-}^{k}\mathbf{,}\nu \right) $, and $%
p_{a+}^{k}\in V\left( \mathbf{A}_{a+}^{k}\mathbf{,}\nu \right) $ it holds that  
$$
p_{a-}^{k}\left( x_{1},\dots,x_{k-1},a,x_{k-1},\dots,x_{n}\right), $$ and 
$$
p_{a+}^{k}\left( x_{1},\dots,x_{k-1},a,x_{k-1},\dots,x_{n}\right) $$ are elements
of the vector space $V\left( \mathbf{A}_{a}^{k}\mathbf{,}\nu \right)$. Therefore the
polynomial $$p:=p_{a-}^{k}\left( x_{1},\dots,x_{k-1},a,x_{k-1},\dots,x_{n}\right)
-p_{a+}^{k}\left( x_{1},\dots,x_{k-1},a,x_{k-1},\dots,x_{n}\right) $$ 
is also an element of $V\left( \mathbf{A}_{a}^{k}\mathbf{,}\nu \right) .$ 
By  (\ref{eqtaged1}) and (\ref{eqtaged2}), for  $~\mathbf{a\in A}_{a}^{k},~%
\mathbf{m\in }\left[ \mathbf{0,}\left( \nu \left( a_{1}\right) ,\dots,\nu
\left( a_{k-1}\right) ,1,\nu \left( a_{k+1}\right) ,\dots,\nu \left(
a_{n}\right) \right) -\mathbf{1}\right] $, we get   $$\partial
_{a_{1}}^{m_{1}}\dots\partial _{a_{k-1}}^{m_{k-1}}\partial
_{a_{k+1}}^{m_{k+1}}\dots\partial _{a_{n}}^{m_{n}}p\left(
x_{1},\dots,x_{k-1},a,x_{k-1},\dots,x_{n}\right) =0.$$ Therefore by Corollary \ref{cor24}, it
follows $p=0.$ 
}
\end{proof}
\textcolor{black}{
Since, after substitution $x_k=a$, the two polynomials $p_{a-}^k$, $p_{a+}^k$ are identical, it follows that their derivatives are also continuous.
Considering the above and Theorem \ref{thm1} about the uniqueness of Hermite interpolation polynomial, it follows that the proposed Hermite spline is regular}. 

\begin{exam} \label{examplesplinecont}
{\rm 
\textcolor{black}{
Let us demonstrate the continuity of the proposed Hermite spline by constructing two Hermite interpolants each one generated by $2 \times 2$ support points on the $\mathbb{R}^2$ plane, while sharing two of these support points. More specifically, let $A_1=\{0,1,2\}$ and $A_2=\{0,1\}$. Then $\mathbf{A}=A_1 \times A_2$. Following the formalism of subsection \ref{contsplines}, we select the support points with $a=1$  (along the $1^{st}$ dimension: $k=1$). Thus, $\mathbf{A}_{1^-}^1=\{0,1\} \times \{0,1\}$  are the support points of Hermite polynomial $p_{1-}^{1}$ and $\mathbf{A}_{a+}^{1} =\{1,2\} \times \{0,1\}$ are the support points of Hermite polynomial $p_{1+}^{1}$. Clearly, the common support points are given by $\mathbf{A}_{a-}^{k}\cap \mathbf{A}_{a+}^{k}=\mathbf{A}_{a}^{k}$, which in this example are just two points: $(1,0)$ and $(1,1)$ both with $x_1$ coordinate equal to 1. The multiplicity function is set for all points as $\nu(\mathbf{a})=3$, for all $\mathbf{a}$.
According to  Theorem \ref{finterpform} the two Hermite interpolants are given by the formulas:
} 
$p_{1-}^{1}(x_1,x_2)=-\left(-5\,{x_{2}}^2+x_{2}+1\right)\,\left({x_{1}}^3-2\,{x_{1}}^2+3\,x_{1}-4\right)$, 

$p_{1+}^{1}(x_1,x_2)=-4.2516\,{x_{1}}^5\,{x_{2}}^5+11.7596\,{x_{1}}^5\,{x_{2}}^4-8.6201\,{x_{1}}^5\,{x_{2}}^3+24.9583\,{x_{1}}^5\,{x_{2}}^2-2.5212\,{x_{1}}^5\,x_{2}-5.2075\,{x_{1}}^5+30.9691\,{x_{1}}^4\,{x_{2}}^5-86.0186\,{x_{1}}^4\,{x_{2}}^4+63.3642\,{x_{1}}^4\,{x_{2}}^3-179.7732\,{x_{1}}^4\,{x_{2}}^2+17.3778\,{x_{1}}^4\,x_{2}+37.5189\,{x_{1}}^4-87.3513\,{x_{1}}^3\,{x_{2}}^5+243.6424\,{x_{1}}^3\,{x_{2}}^4-180.4411\,{x_{1}}^3\,{x_{2}}^3+494.6803\,{x_{1}}^3\,{x_{2}}^2-45.0001\,{x_{1}}^3\,x_{2}-103.1784\,{x_{1}}^3+118.7548\,{x_{1}}^2\,{x_{2}}^5-332.4123\,{x_{1}}^2\,{x_{2}}^4+247.3391\,{x_{1}}^2\,{x_{2}}^3-649.9843\,{x_{1}}^2\,{x_{2}}^2+54.9454\,{x_{1}}^2\,x_{2}+135.4968\,{x_{1}}^2-78.0744\,x_{1}\,{x_{2}}^5+219.1735\,x_{1}\,{x_{2}}^4-163.7111\,x_{1}\,{x_{2}}^3+420.2292\,x_{1}\,{x_{2}}^2-33.7958\,x_{1}\,x_{2}-87.4964\,x_{1}+19.9533\,{x_{2}}^5-56.1446\,{x_{2}}^4+42.0691\,{x_{2}}^3-120.1103\,{x_{2}}^2+10.9939\,x_{2}+24.8667$

\textcolor{black}{
Figure \ref{figsplincont2d} shows $p_{1-}^{1}(x_1,x_2)$ and $p_{1+}^{1}(x_1,x_2)$ in blue and green color respectively.
By substituting $x_1 = 1$ in $p_{1-}^{1}, p_{1+}^{1}$ we get 
\begin{eqnarray*}
    p_{1-}^{1} = -10x_2^2+ 2 x_2 +2, \;\;\;p_{1+}^{1}(x_2)= -10x_2^2+ 2 x_2 +2,
\end{eqnarray*}
which confirms the continuity of the two polynomials and their derivatives.
}

\begin{figure}[!h]
     \centering
         \includegraphics[width=9cm]{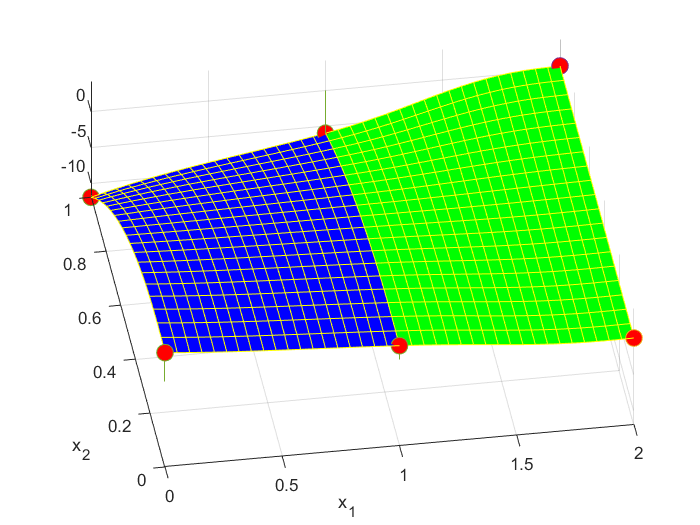}
         \caption{The two Hermite interpolants of example \ref{examplesplinecont} with the support points shown as red circles. The continuity of the Hermite spline along the points with $x_1=1$ is visually confirmed.}
         \label{figsplincont2d}
\end{figure}
}
\end{exam}

 %%%%%%%%%%%%%%%%%%%%%%%%%%%%%%%%%

\textcolor{black}{In the example below we provide a univariate example to further facilitate the continuity demonstration. More specifically we generate a univariate Hermite spline over $7$ support points for each Hermite interpolant, using $4$ local support points, for different values of the multiplicity function $v(a)=1,2,3$. The continuity of the interpolants and their derivative can be visually confirmed.}

\begin{exam} \label{exxxx}
{\rm
\textcolor{black}{
Consider the grid $A=\{0,\dots ,6\}$ and the Hermite polynomial $f(x)$ which satisfies $\partial _{{a}}^{{m}}f=t_{{a}}^{{m}}, 
{a\in A}$, ${m\in }\left[ {0,}\nu \left( {a}%
\right) -1\right]$ with known values given in Table \ref{tablyexxx}.
For any $x \in [0,6]$ we derive the spline Hermite polynomial, defined in $\Big[\lfloor x \rfloor, \lceil x \rceil \Big]$, using $|B|=4$ number of local support points, with $B=\Big[\lfloor x \rfloor -1, \lceil x \rceil +1 \Big]$
%Let $g=(x-2)^8 - (x+3)^4 +(x-7)^3+7 x^2+8$,
for different values of the multiplicity function $\nu(a) = \{1,2,3\} $ and for every point $a$ on the grid $ A$. In Fig. \ref{fisubfigxx}(a) the spline Hermite polynomials are shown in different colors for each domain, whereas the style of the line indicates the value of $\nu(a)$, as follows: dotted line: $\nu(a)=1$, dashed line: $\nu(a)=2$ and continuous line: $\nu(a)=3$. 
The 1st, 2nd, and 3rd derivatives of the spline Hermite polynomials are shown in Fig.  \ref{fisubfigxx}(b), (c), and (d), respectively. As was expected, the Hermite spline is continuous ($C^0$) for all values of $\nu(a)>0$, the 1st derivative of the Hermite spline is continuous ($C^1$) for $\nu(a)=\{2,3\}$, the 2nd derivative is continuous ($C^2$) only for $\nu(a)=3$.}

\begin{table}[!h]
\begin{center}
\caption{The conditions that must be satisfied by the Hermite polynomial in Example  \ref{exxxx}}\label{tablyexxx}
\begin{tabular}{|l ccc|}
\hline
$\mathbf{a}$ & $t^0_{{a}}$ & $t^1_{{a}}$ & $t^2_{{a}}$\\ \hline
0 &  -160 & -985 & 3448\\ \hline
1  &  -456  &  -142 & -158\\ \hline
2 &  -714 &  -397 & -316 \\ \hline
3 &  -1288 &  -766 & -386 \\ \hline
4 & -2052 &   -265 & 2992\\ \hline
5   & 2640 &  15530 & 40058\\ \hline  
6  &  59234 & 128243 & 228412\\ \hline
\hline
\end{tabular}
\end{center} 
\end{table}

\begin{figure}[h!]
\centering
\subfigure[]{{\includegraphics[width=0.5\textwidth]{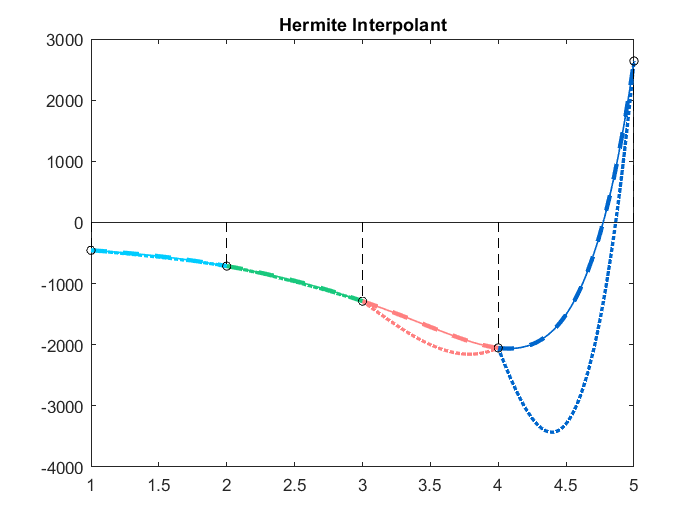}}}\hfill
\subfigure[]{{\includegraphics[width=0.5\textwidth]{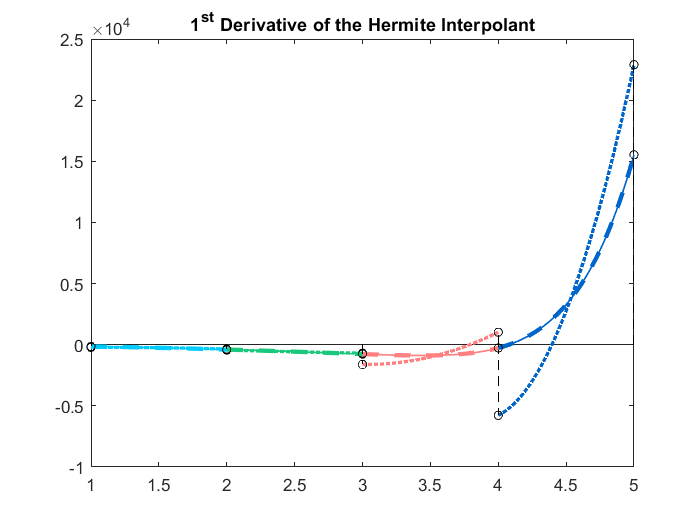}}}
\subfigure[]{{\includegraphics[width=0.5\textwidth]{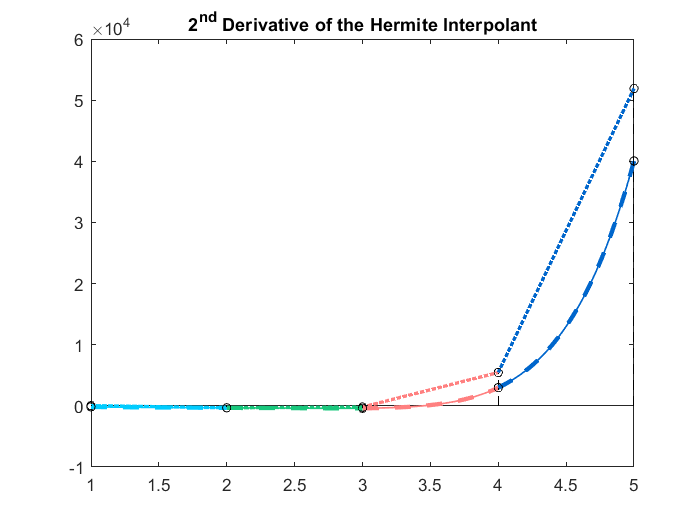}}}\hfill
\subfigure[]{{\includegraphics[width=0.5\textwidth]{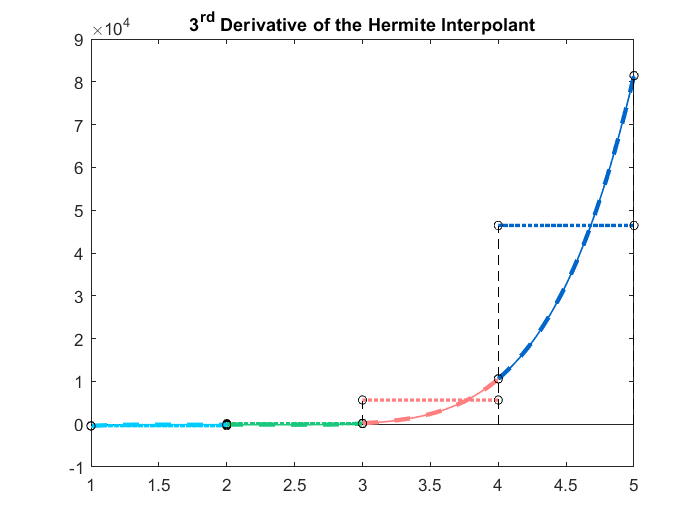}}}
\caption{\textcolor{black}{The continuity of the interpolants of $f(x)$, with multiplicity functions $\nu(a) = \{1,2,3\} $, for every point $a$ on the grid $ A$,  and its derivatives up to order 3 of Example \ref{exxxx}. The different colors indicate the different spline interpolants for each interval. Dotted, dashed and continuous line indicates the interpolant created with $v(a)=1$, $v(a)=2$ and $v(a)=3$, respectively.  (a) the $C^0$ continuity of the spline interpolants $f(x)$ (b) $C^1$ continuity of the spline interpolants of $f(x)$ (c) the $C^2$ continuity of the spline interpolants of $f(x)$ (d) the $C^3$ continuity of the spline interpolants of $f(x)$}.}
\label{fisubfigxx}
\end{figure}

}
\end{exam}

%\newpage
%\subsection{Continuity of the Splines}
%\label{contsplines}
%\input{splinecontinunity}

\newpage
\section{Examples}
\label{examm}
%%%%%%%%%%%%%%%%%%%%%%%%%%%%%%%%%%%%%%%%%%
%%%%%%%%%%%%%%%%%%%%%%%%%%%%%%%%%%%%%%%%%%

\begin{exam}
{\rm
Consider the grid $\mathbf{A}=A_1 \times A_2$, where $A_1 = \{0,1\},\; A_2 = \{0,1\}$, with multiplicity function $v(a) = (v_1(a),\;v_2(a)) = (2,2)$, for every point on the grid, $a \in \mathbb{A}$. Given the values of the interpolating polynomial $f$, $t_{\mathbf{a}%
}^{\mathbf{k}},~\mathbf{a\in A},~\mathbf{k\in }\left[ \mathbf{0,}\nu \left( 
\mathbf{a}\right) -\mathbf{1}\right] $, namely $t_{(0,0)}^{(0,0)},t_{(1,0)}^{(0,0)},t_{(0,1)}^{(0,0)},t_{(1,1)}^{(0,0)}$ at each one of the support points, as well as the values of the partial derivatives up to the 1st degree: $t_{(0,0)}^{(0,1)},t_{(1,0)}^{(0,1)},t_{(0,1)}^{(0,1)},t_{(1,1)}^{(0,1)}$, $t_{(0,0)}^{(1,0)},t_{(1,0)}^{(1,0)},t_{(0,1)}^{(1,0)},t_{(1,1)}^{(1,0)}$, $t_{(0,0)}^{(1,1)},t_{(1,0)}^{(1,1)},t_{(0,1)}^{(1,1)},t_{(1,1)}^{(1,1)}$, the following conditions must hold: 
$$\partial _{\mathbf{a}}^{\mathbf{k}}f=t_{\mathbf{a}%
}^{\mathbf{k}},~\mathbf{a\in A},~\mathbf{k\in }\left[ \mathbf{0,}\nu \left( 
\mathbf{a}\right) -\mathbf{1}\right].$$

We compute the $H_{\left( i,a\right) }\left( x_{i}\right) $ quantities given by (\ref{Hmikra}) for $i\in \left\{ 1,2\right\} $ and $a\in A_{i}$:
\begin{eqnarray*}
&H_{(1,0)} =x_1^2 - 2x_1+1&,H_{(2,0)} = x_2^2 -2x_2 + 1\\ &H_{(1,1)}=x_1^2&, H_{(2,1)} = x_2^2
\end{eqnarray*}

Furthermore, we compute the $H_{(\mathbf{a},\mathbf{k})}$ quantities  given by (\ref{Hmegala}) for every $\mathbf{a} \in \mathbf{A}$ and $\mathbf{k} \in [\mathbf{0},v(\mathbf{a})-\mathbf{1}]$. The numbering of the elements of $[\mathbf{0},v(\mathbf{a})-\mathbf{1}]$ is given by the degree reverse lexicographic order. More specifically, $$
1_{\mathbf{a}} =\left( 0,0\right) \prec 2_{\mathbf{a}}=\left(
0,1\right) \prec 3_{\mathbf{a}}=\left(1,0\right)\prec  4_{\mathbf{a}} =\left( 1,1\right).$$
\begin{eqnarray*}
&H_{(0,0),(0,0)}= (x_1^2 - 2x_1+1)(x_2^2 - 2x_2+1)&, H_{(0,1),(0,0)}= (x_1^2 - 2x_1+1)x_2^2\\ 
&H_{(0,0),(0,1)}= (x_1^2 - 2x_1+1)(x_2^3 - 2x_2^2+x_2)&,
H_{(0,1),(0,1)}= (x_1^2 - 2x_1+1)(x_2^3-x^2)\\
&H_{(0,0),(1,0)}= (x_1^3 - 2x_1^2+x_1)(x_2^2 - 2x_2+1)&,
H_{(0,1),(1,0)}= (x_1^3 - 2x_1^2+x_1)x_2^2\\
&H_{(0,0),(1,1)}= (x_1^3 - 2x_1^2+x_1)(x_2^3 - 2x_2^2+x_2)&,
H_{(0,1),(1,1)}= (x_1^3 - 2x_1^2+x_1)(x_2^3-x^2)\\
\\
&H_{(1,0),(0,0)}= x_1^2(x_2^2 - 2x_2+1)&,
H_{(1,1),(0,0)}= x_1^2x_2^2\\
&H_{(1,0),(0,1)}= x_1^2(x_2^3 - 2x_2^2+x_2)&,
H_{(1,1),(0, 1)}=x_1^2(x_2^3-x_2^2)\\
&H_{(1,0),(1,0)}=(x_1^3-x_1^2)(x_2^2 - 2x_2+1)&,
H_{(1,1),(1,0)}= (x_1^3-x_1^2)x_2^2\\
&H_{(1,0),(1,1)}= (x_1^3-x_1^2)(x_2^3 - 2x_2^2+x_2)&,
H_{(1,1),(1,1)}= (x_1^3-x_1^2)(x_2^3-x_2^2)\\
\end{eqnarray*}

Thus, $H_{{(0,0)}}=\left(H_{(0,0),(0,0)},\; H_{(0,0),(0,1)},\; H_{(0,0),(1,0)},\; H_{(0,0),(1,1)} \right)^T$. Similarly, for $H_{{(0,1)}},\; H_{{(1,0)}}$ and $H_{{(1,1)}}$.

For  each point on the grid, we create the lower triangular matrix, $L_{\mathbf{a}}$:
\begin{eqnarray*}
    \Lambda_{(0,0)}^{-1} =
    \left(  
    \begin{array}{cccc}
    1 &  0 & 0 & 0 \\
    -2 & 1 & 0 & 0 \\
    -2 & 0 & 1 & 0 \\
    4 & -2 & -2 & 1
    \end{array}   
    \right)&,&
        \Lambda_{(0,1)}^{-1} =
    \left(  
    \begin{array}{cccc}
    1 &  0 & 0 & 0 \\
    2 & 1 & 0 & 0 \\
    -2 & 0 & 1 & 0 \\
    -4 & -2 & 2 & 1
    \end{array}   
    \right),\\
    \\
    \Lambda_{(1,0)}^{-1} =
    \left(  
    \begin{array}{cccc}
    1 &  0 & 0 & 0 \\
    -2 & 1 & 0 & 0 \\
    2 & 0 & 1 & 0 \\
    -4 & 2 & -2 & 1
    \end{array}   
    \right)&,&
        \Lambda_{(1,1)}^{-1} =
    \left(  
    \begin{array}{cccc}
    1 &  0 & 0 & 0 \\
    2 & 1 & 0 & 0 \\
    2 & 0 & 1 & 0 \\
    4 & 2 & 2 & 1
    \end{array}   
    \right),
\end{eqnarray*}

The interpolating polynomial obtained by using the formula of Theorem \ref{finterpform} is
\begin{eqnarray}\label{exampinterpol}
     f &=& \left(  \Lambda_{(0,0)}^{-1} \cdot
    \left(  
    \begin{array}{c}
    t_{(0,0)}^{(0,0)} \\
    t_{(0,0)}^{(0,1)}  \\
    t_{(0,0)}^{(1,0)}  \\
    t_{(0,0)}^{(1,1)} 
    \end{array}   
    \right)\right)^T
     \left(  
    \begin{array}{c}
    H_{(0,0),(0,0)} \\
    H_{(0,0),(0,1)}  \\
    H_{(0,0),(1,0)} \\
    H_{(0,0),(1,1)} 
    \end{array}   
    \right) 
    +
    \left(  \Lambda_{(0,1)}^{-1} \cdot
    \left(  
    \begin{array}{c}
    t_{(0,1)}^{(0,0)} \\
    t_{(0,1)}^{(0,1)}  \\
    t_{(0,1)}^{(1,0)}  \\
    t_{(0,1)}^{(1,1)} 
    \end{array}   
    \right)\right)^T
     \left(  
    \begin{array}{c}
    H_{(0,1),(0,0)} \\
    H_{(0,1),(0,1)}  \\
    H_{(0,1),(1,0)} \\
    H_{(0,1),(1,1)} 
    \end{array}   
    \right) \nonumber\\
    &+&
    \left(  \Lambda_{(1,0)}^{-1} \cdot
    \left(  
    \begin{array}{c}
    t_{(1,0)}^{(0,0)} \\
    t_{(1,0)}^{(0,1)}  \\
    t_{(1,0)}^{(1,0)}  \\
    t_{(1,0)}^{(1,1)} 
    \end{array}   
    \right)\right)^T
     \left(  
    \begin{array}{c}
    H_{(1,0),(0,0)} \\
    H_{(1,0),(0,1)}  \\
    H_{(1,0),(1,0)} \\
    H_{(1,0),(1,1)} 
    \end{array}   
    \right) 
    +
    \left(  \Lambda_{(1,1)}^{-1} \cdot
    \left(  
    \begin{array}{c}
    t_{(1,1)}^{(0,0)} \\
    t_{(1,1)}^{(0,1)}  \\
    t_{(1,1)}^{(1,0)}  \\
    t_{(1,1)}^{(1,1)} 
    \end{array}   
    \right)\right)^T
     \left(  
    \begin{array}{c}
    H_{(1,1),(0,0)} \\
    H_{(1,1),(0,1)}  \\
    H_{(1,1),(1,0)} \\
    H_{(1,1),(1,1)} 
    \end{array}   
    \right) .
\end{eqnarray}
}
\end{exam}
%%%%%%%%%%%%%%%%%%%%%%%%%%%%%%%%%%%%%%%%%%
%%%%%%%%%%%%%%%%%%%%%%%%%%%%%%%%%%%%%%%%%%
\begin{exam}
{\rm
Let us expand the previous example by using the proposed Hermite polynomial in order to interpolate a given function $g(x_1,x_2) = e^{x_1+x_2}$. To this end, it is sufficient to 
set the values of $t_{ij}^{(km)}$ equal to $\partial_{ij}^{(km)}g$. More specifically:

\begin{equation*}
    \left(  
    \begin{array}{c}
    t_{(0,0)}^{(0,0)} \\
    t_{(0,0)}^{(0,1)}  \\
    t_{(0,0)}^{(1,0)}  \\
    t_{(0,0)}^{(1,1)} 
    \end{array}\right) =
        \left(  
    \begin{array}{c}
    1 \\
    1\\
    1\\
    1
    \end{array}\right),
        \left(  
    \begin{array}{c}
    t_{(0,1)}^{(0,0)} \\
    t_{(0,1)}^{(0,1)}  \\
    t_{(0,1)}^{(1,0)}  \\
    t_{(0,1)}^{(1,1)} 
    \end{array}\right)
    =
    \left(  
    \begin{array}{c}
    t_{(1,0)}^{(0,0)} \\
    t_{(1,0)}^{(0,1)}  \\
    t_{(1,0)}^{(1,0)}  \\
    t_{(1,0)}^{(1,1)} 
    \end{array}\right)
    =
        \left(  
    \begin{array}{c}
    e \\
    e\\
    e\\
    e
    \end{array}\right),
        \left(  
    \begin{array}{c}
    t_{(1,1)}^{(0,0)} \\
    t_{(1,1)}^{(0,1)}  \\
    t_{(1,1)}^{(1,0)}  \\
    t_{(1,1)}^{(1,1)} 
    \end{array}\right) =
        \left(  
    \begin{array}{c}
    e^2 \\
    e^2\\
    e^2\\
    e^2
    \end{array}\right).
\end{equation*}
Substituting the values of $t_{ij}^{(km)}$ in (\ref{exampinterpol}) of the previous example and performing some algebra, yields the following Hermite polynomial. 
\begin{equation*}
    f(x_1,x_2)= (1 - x_1(-1 + (5 + e(-2 + x_1) - 3x_1)x_1))(1 + x_2 + (-5 + 2e x_2^2 - (-3 + e)x_2^3)
\end{equation*}
This polynomial is identical with the one obtained by Eq. (\ref{lemrel1}) from Lemma \ref{lemErr}. 

The calculated RMSE between the original function and the interpolated polynomial, in a grid of $11 \times 11$ equally spaced points in $[0,1] \times [0,1]$ is $0.0085$, which suggests a very accurate interpolation.
}
\end{exam}
%%%%%%%%%%%%%%%%%%%%%%%%%%%%%%%%%%%%%%%%%%
%%%%%%%%%%%%%%%%%%%%%%%%%%%%%%%%%%%%%%%%%%

\begin{exam}\label{examp2d}
{\rm
Consider the function $g(x_1,x_2) = e^{-(x_1 - 3)^2-(x_2 -3)^2} + e^{\frac{-(x_1 - 4)^2-(x_2-4)^2}{5}}$ and the grid of support points $\mathbf{A}=A_1 \times A_2$, where $A_1 = \{0 ,1,2,3,4,5\},\; A_2 = \{0 ,1,2,3,4,5\}$, with multiplicity function (maximum order of derivation -1) $v(\mathbf{a}) = (v_1(a),\;v_2(a))$, for every point on the grid, $a \in \mathbb{A}$. Thus in total  $6\times 6$ support points have been defined.
The 3d plot of the function, $g(x_1,x_2)$ is depicted in Figure \ref{fig1}. In addition, the difference between $g(x_1,x_2)$ and the interpolated function using the proposed Hermit with maximum order of derivative $v(\mathbf{a}) = (v_1(a),\;v_2(a)) = (2,2)$ and the cubic spline interpolation (which was the best performer of the methods under comparison) are shown in Figure   \ref{fisubfig}(a) and (b) respectively. Note the the color scale is identical in all difference images.
\begin{figure}[!h]
     \centering
         \includegraphics[width=9cm]{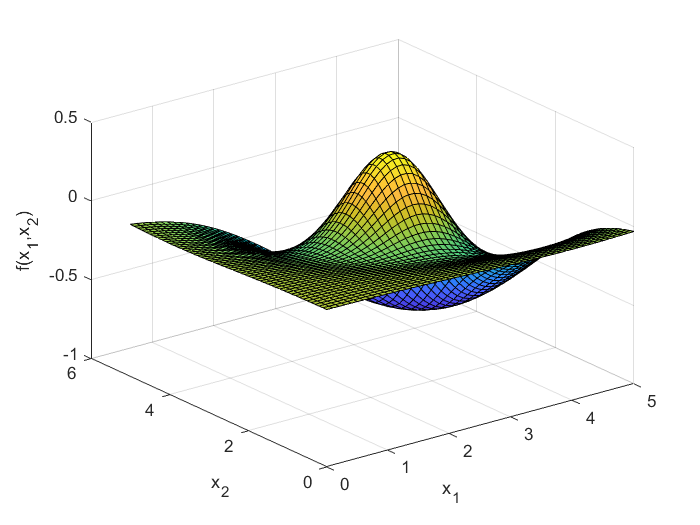}
         \caption{Function $g(x_1,x_2) = e^{-(x_1 - 3)^2-(x_2 -3)^2} + e^{\frac{-(x_1 - 4)^2-(x_2-4)^2}{5}}$ of Example \ref{examp2d}}
         \label{fig1}
\end{figure}

\begin{figure}[h!]
\centering
\subfigure[]{{\includegraphics[width=0.5\textwidth]{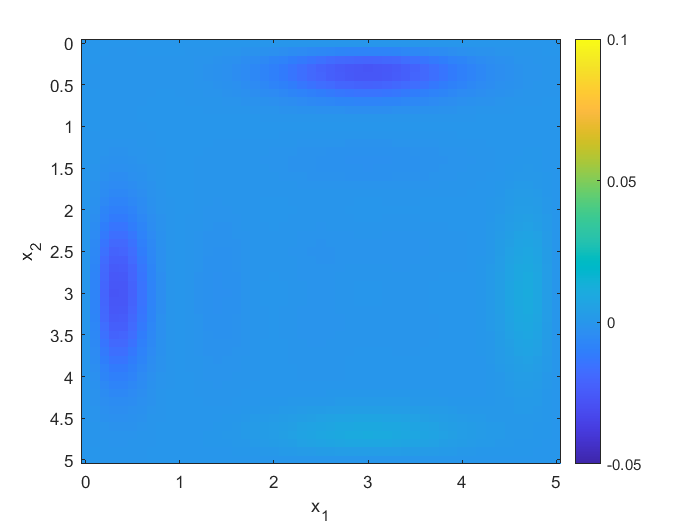}}}\hfill
\subfigure[]{{\includegraphics[width=0.5\textwidth]{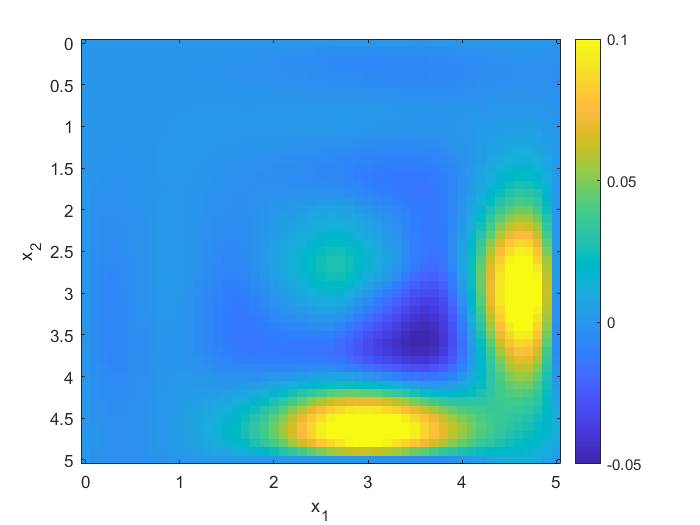}}}
\caption{The difference between $g(x_1,x_2)$ of Example \ref{examp2d} and the interpolated function using (a) the proposed Hermit with multiplicity function $v(a) = (v_1(a),\;v_2(a)) = (2,2)$ and (b) the cubic spline interpolation. Note the the color scale is identical in all difference images.}
\label{fisubfig}
\end{figure}

In this example, we utilize the values of $g$ and its partial derivatives up to maximum order (m.o.) $v(\mathbf{a})-1$ at each support point $\mathbf{a}$, to interpolate the function at $51 \times 51$ equally distributed points in $[0,5] \times [0,5]$.\textcolor{black}{ Please note that the proposed Hermite is a single polynomial for the entire domain, rather than piece-wise polynomials. The cubic spline $C^2$, linear spline $C^0$, cubic spline $C^1$, and the Hermit modified Akima $C^1$ ("Makima",\cite{practspline}) interpolation methods are also applied at the same points, using the Matlab implementation.} The achieved RMSE by each method is shown in Table \ref{tably1}, calculated at the aforementioned points. These interpolation methods have been used for comparison purposes in the remaining examples.
\begin{table}[h!]
\begin{center}
\caption{The achieved RMSE for the proposed Hermite polynomial and other interpolation methods, in the case of Example \ref{examp2d}}\label{tably1}
\begin{tabular}{|l cc|}
\hline
\textbf{Method} & \textbf{RMSE} &\\ \hline
Proposed Hermite $v(\mathbf{a}) = (2,2)$ & 0.0054  & \\ \hline
Proposed Hermite $v(\mathbf{a})= (3,3)$ & 0.0002 & \\ \hline
\textcolor{black}{Linear Spline} &  0.0554 &  \\ \hline
\textcolor{black}{Cubic Spline $C^1$}   & 0.0237 &  \\ \hline
Makima  & 0.0255 & \\ \hline
\textcolor{black}{Cubic Spline $C^2$}   & 0.0223 & \\ \hline
\end{tabular}
\end{center}
\end{table}

}
\end{exam}

%%%%%%%%%%%%%%%%%%%%%%%%%%%%%%%%%%%%%%%%%%
%%%%%%%%%%%%%%%%%%%%%%%%%%%%%%%%%%%%%%%%%%
\begin{exam}\label{examp3d}
{\rm
Consider the function $g(x_1,x_2,x_3) = e^{\frac{-(x_1 - 3)^2-(x_2 -1)^2-(x_3 -1.5)^2)}{3}}- e^{\frac{-(x_1 - 0.5)^2-(x_2 -2)^2-(x_3 -1)^2)}{5}}$, and the grid of support points $\mathbf{A}=A_1 \times A_2 \times A_3$, where $A_1 = \{0 ,1,2,3\},\; A_2 = \{0 ,1,2,3,4\}$ and $A_3 = \{0,1,2\}$, with multiplicity function $v(\mathbf{a}) = (v_1(\mathbf{a}),\;v_2(\mathbf{a}),\;v_3(\mathbf{a})) = (1,1,1) , (2,2,2)$ and $(3,3,3)$, for every point on the grid, $\mathbf{a} \in \mathbf{A}$. Thus, in total, we utilize $4 \times 5 \times 3$ support points to derive \textcolor{black}{ a single Hermite interpolating polynomial for the entire domain}. 

In this example we use the values of $g$ and its partial derivatives up to maximum order,  $v(\mathbf{a})-1$, at each support point $\mathbf{a}$, to interpolate the function at $13 \times 17 \times 9$ equally distributed points in $[0,3] \times [0,4] \times [0,2]$. The cubic spline, linear, cubic polynomial, and the Hermit Makima interpolation methods are also applied at the same points. The achieved RMSE by each method is shown in Table \ref{table2}, calculated at the aforementioned points.

\begin{table}[!h]
\begin{center}
\caption{The achieved RMSE for the proposed Hermite polynomial and other interpolation methods, in the case of Example \ref{examp3d}}\label{table2}
\begin{tabular}{|l cc|}
\hline
\textbf{Method} & \textbf{RMSE} &\\ \hline
Proposed Hermite $v(\mathbf{a}) = (1,1,1)$ & 0.0152  & \\ \hline
Proposed Hermite $v(\mathbf{a}) = (2,2,2)$ &  0.0001 & \\ \hline
Proposed Hermite $v(\mathbf{a}) = (3,3,3)$ & 1.2776e-06 & \\ \hline
\textcolor{black}{Linear Spline} &  0.0528 &  \\ \hline
\textcolor{black}{Cubic Spline $C^1$}   & 0.0187 &  \\ \hline
Makima  & 0.0209 & \\ \hline
\textcolor{black}{Cubic Spline $C^2$}   & 0.0152 & \\ \hline
\end{tabular}

\end{center} 
\end{table}
}
\end{exam}

%%%%%%%%%%%%%%%%%%%%%%%%%%%%%%%%%%%%%%%%%%
%%%%%%%%%%%%%%%%%%%%%%%%%%%%%%%%%%%%%%%%%%

\begin{exam}
{\rm
    Consider the grid $\mathbf{A}=A_1 \times A_2$, where $A_1 = \{0,1\},\; A_2 = \{0,1\}$, with multiplicity function $v(a) = (v_1(a),\;v_2(a)) = (1,1)$, for every point on the grid, $a \in \mathbf{A}$. The given values of the support points are $c_{00}, \; c_{01},\;c_{10},\;c_{11}$, such that $t_{(i,j)}^{(0,0)}=c_{ij}, i,j = \{0,1\}$.

    The expressions $H_{\left( i,a\right) }\left( x_{i}\right) $ in (\ref{Hmikra}) are obtained for $i\in \left\{ 1,2\right\} $ and $a\in A_{i}$ as following:  $H_{(1,0)} =1-x_1,\; H_{(1,1)} =x_1,\; H_{(2,0)} =1-x_2$, and $H_{(2,1)} =x_2$.

    Furthermore, we derive $H_{(\mathbf{a},\mathbf{k})}$ given by (\ref{Hmegala}) for every $\mathbf{a} \in \mathbf{A}$ and $\mathbf{k}=(0,0)$: 
    \begin{eqnarray*}
     H_{(0,0),(0,0)} &=& (1-x_1)(1-x_2) ,\\ H_{(0,1),(0,0)}&=& (1-x_1)x_2,\\ H_{(1,0),(0,0)}&=& x_1(1-x_2), \\ H_{(1,1),(0,0)} &=& x_1x_2.
    \end{eqnarray*}

    Note that for every $\mathbf{a} \in \mathbf{A}$, $\Lambda_{\mathbf{a}}=1$. Thus, the interpolating polynomial is given by the following expression:
    \begin{equation}\label{bil}
        f(x_1,x_2) = c_{00} (1-x_1)(1-x_2)+c_{01}(1-x_1)x_2 + c_{10} x_1(1-x_2)  + c_{11}x_1x_2
    \end{equation}
    
    Observe that the above expression is identical to the bilinear interpolating polynomial for the unit step grid.
    
    %\begin{comment}
    %\begin{tikzpicture}
    %\draw[help lines, %color=gray!30, dashed] %(0,0) grid (2,2);
    %\draw[->,ultra thick] %(-0.1,0)--(2,0) %node[right]{$x_1$};
    %\draw[->,ultra thick] (0,-0.1)--(0,2) node[above]{$x_2$};
    %\end{tikzpicture}
    %\end{comment} 
}
\end{exam}

%%%%%%%%%%%%%%%%%%%%%%%%%%%%%%%%%%%%%%%%%%
%%%%%%%%%%%%%%%%%%%%%%%%%%%%%%%%%%%%%%%%%%

\begin{exam}
    Consider the grid $\mathbf{A}=A_1 \times A_2 \times A_3$, where $A_1 = \{0,1\},\; A_2 = \{0,1\}, \; A_3 = \{0,1\}$, with multiplicity function $v(a) = (v_1(a),\;v_2(a),\;v_3(a)) = (1,1,1)$, for every point on the grid, $a \in \mathbf{A}$. The given values of the support points are $c_{000},c_{001},c_{010},c_{100},c_{011},c_{101},c_{110},c_{111}$, such that $t_{(i,j,k)}^{(0,0,0)}=c_{ijk}, i,j,k = \{0,1\}$.

    We derive $H_{\left( i,a\right) }\left( x_{i}\right) $  given by (\ref{Hmikra}) for $i\in \left\{ 1,2\right\} $ and $a\in A_{i}$:  $H_{(1,0)} =1-x_1,\; H_{(1,1)} =x_1,\; H_{(2,0)} =1-x_2$,$H_{(2,1)} =x_2, H_{(3,0)} =1-x_3,$ and $\; H_{(1,1)} =x_3$.
    Furthermore, we obtain $H_{(\mathbf{a},\mathbf{k})}$   given by (\ref{Hmegala}) for every $\mathbf{a} \in \mathbf{A}$ and $\mathbf{k}=(0,0,0)$: 
\begin{eqnarray*}
&H_{(0,0,0),(0,0,0)}= (1-x_1)(1-x_2)(1-x_3)&, H_{(0,0,1),(0,0,0)}=  (1-x_1)(1-x_2)x_3\\ 
&H_{(0,1,0),(0,0,0)}=  (1-x_1)x_2(1-x_3)&,
H_{(1,0,0),(0,0,0)}= x_1(1-x_2)(1-x_3)\\
&H_{(0,1,1),(0,0,0)}=(1-x_1)x_2x_3&,
H_{(1,0,1),(0,0,0)}= x_1(1-x_2)x_3\\
&H_{(1,1,0),(0,0,0)}= x_1 x _2 (1-x_3)&,
H_{(1,1,1),(0,0,0)}= x_1 x_2 x_3
\end{eqnarray*}

%For every $\mathbf{a} \in \mathbf{A}$, we substitute $\Lambda_{\mathbf{a}}=1$ from (\ref{La}), quantities as and the expression of from (4.12) into the inequality (4.7) to deduce (\ref{mainform})
Note that for every $\mathbf{a} \in \mathbf{A}$, $\Lambda_{\mathbf{a}}=1$. Thus, the interpolating polynomial is given by the following expression:
\begin{eqnarray*}
\label{tril}
        f(x_1,x_2,x_3) &=& c_{000}(1-x_1)(1-x_2)(1-x_3)+c_{001}(1-x_1)x_2(1-x_3)  \\&& +  c_{010}(1-x_1)x_2(1-x_3) 
        + c_{100}x_1(1-x_2)(1-x_3) + c_{011}(1-x_1)x_2x_3\\&& 
        + c_{101} x_1 x_2 (1-x_3) 
        +c _{110}x_1x_2(1-x_3) +c_{111}x_1x_2x_3
    \end{eqnarray*}
    
    Observe that the above expression is identical to the trilinear interpolating polynomial  for the unit step grid.
    
\end{exam}

%%%%%%%%%%%%%%%%%%%%%%%%%%%%%%%%%%%%%%%%%%
%%%%%%%%%%%%%%%%%%%%%%%%%%%%%%%%%%%%%%%%%%
\begin{exam}\label{gridfun3}
{\rm

Consider the function %$\displaystyle 
$ g(x_1,x_2,x_3) =  x_1 \sin(x_2)+ \frac{x_2 \sin(x_1)}{10}-x_1 \sin\left(\frac{x_2 x_3}{4}\right)$, and the grid of support points $\mathbf{A}=A_1 \times A_2 \times A_3$, where $A_1 =A_2=A_3= \{-7,-6,\dots,7\}$, with multiplicity function $v(\mathbf{a}) = (v_1(\mathbf{a}),\;v_2(\mathbf{a}),\;v_3(\mathbf{a})) = (1,1,1) , (2,2,2)$, for every point on the grid, $\mathbf{a} \in \mathbf{A}$. Thus, in total we utilize $15 \times 15 \times 15$ support points to derive the Hermite interpolating polynomial. 

\textcolor{black}{
In this example we use the values of $g$ and its partial derivatives up to maximum order,  $v(\mathbf{a})-1$, at each support point $\mathbf{a}$, to interpolate the function at $57 \times 57 \times 57$ equally distributed points in $[-7,7] \times [-7,7] \times [-7,7]$.  The achieved RMSE by each method is shown in Table \ref{table4}, calculated at the aforementioned points.}

In addition, we implement the Hermit interpolation via splines, where we use a smaller grid. Namely, for the current example, the grid of spline support points $\mathbf{B}=B_1 \times B_2 \times B_3$, where $B_1 =B_2=B_3= \{0,1,2\}$, with multiplicity function $v(\mathbf{b}) = (v_1(\mathbf{b}),\;v_2(\mathbf{b}),\;v_3(\mathbf{b})) = (2,2,2)$, for every point on the grid, $b \in \mathbf{B}$.

\begin{table}[!h]
\begin{center}
\caption{The achieved RMSE for the proposed Hermite polynomial and other interpolation methods, in the case of Example \ref{examp3d}}\label{table4}
\begin{tabular}{|l cc|}
\hline
\textbf{Method} & \textbf{RMSE} &\\ \hline
%Proposed Hermite $v(\mathbf{a}) = (1,1,1)$ &   & \\ \hline
Proposed Hermite $v(\mathbf{a}) = (2,2,2)$ &   0.0015 & \\ \hline
\textcolor{black}{Linear Spline} & 1.3318 &  \\ \hline
\textcolor{black}{Cubic Spline $C^1$}    &  0.6699 &  \\ \hline
Makima  &  0.7748 & \\ \hline
\textcolor{black}{Cubic Spline $C^2$}   &    0.7461 & \\ \hline
\end{tabular}

\end{center} 
\end{table}

\begin{figure}[!h]
     \centering
         \includegraphics[width=16cm]{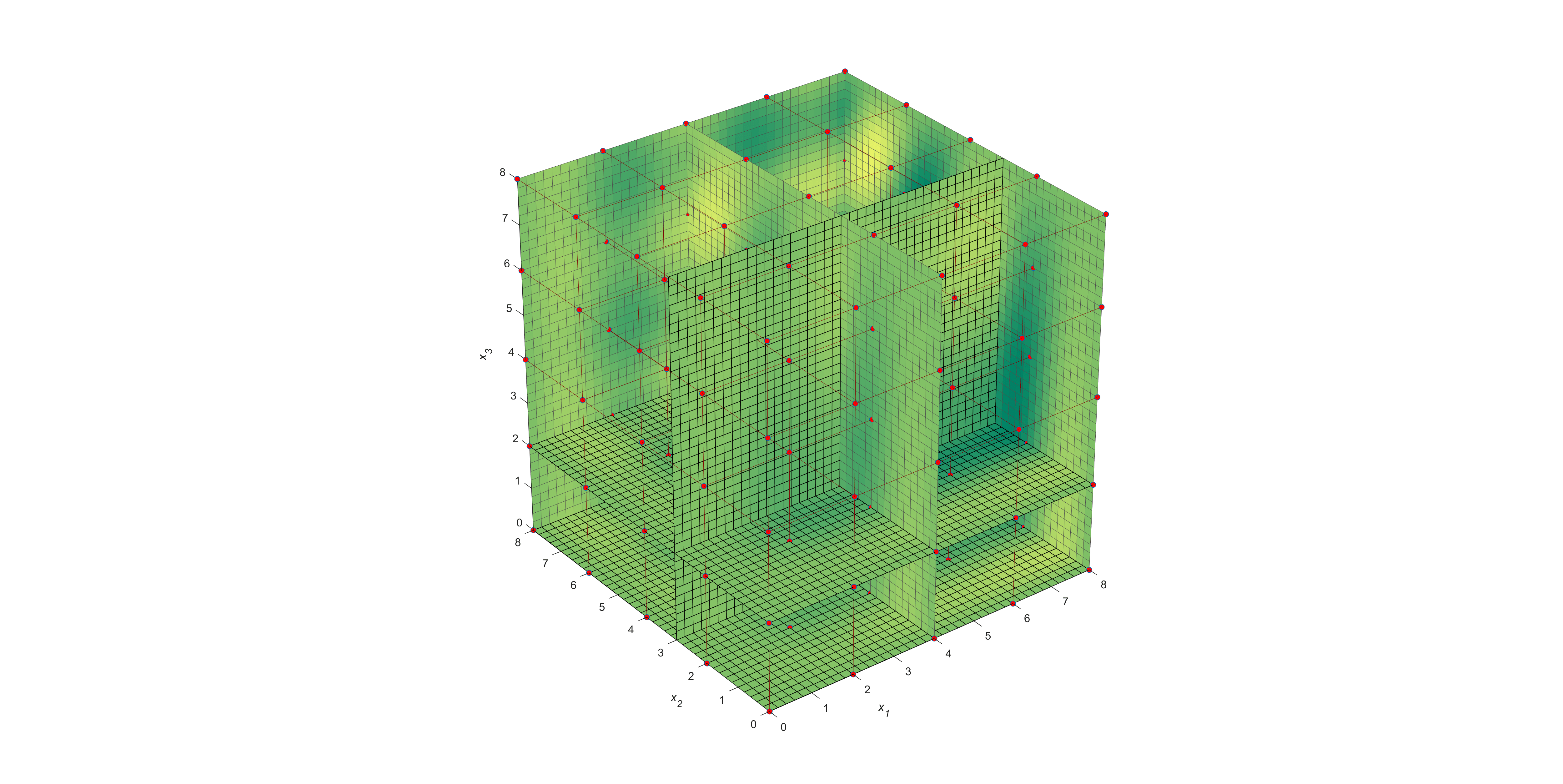}%{3d_1_250.png}
         \caption{Function $ g(x_1,x_2,x_3) =  x_1 \sin(x_2)+ \frac{x_2 \sin(x_1)}{10}-x_1 \sin\left(\frac{x_2 x_3}{4}\right)$ of Example \ref{gridfun3}}
         \label{figgr}
\end{figure}

} 
\end{exam}

%%%%%%%%%%%%%%%%%%%%%%%%%%%%%%%%%%%%%%%%%%%
%%%%%%%%%%%%%%%%%%%%%%%%%%%%%%%%%%%%%%%%%%%
%\newpage
\begin{exam}\label{examp8}
{\rm
Consider the function of Example \ref{gridfun3}, $ g(x_1,x_2,x_3) =  x_1 \sin(x_2)+ \frac{x_2 \sin(x_1)}{10}-x_1 \sin\left(\frac{x_2 x_3}{4}\right)$ and the support points arranged on the same 3D grid. The task is to interpolate the values of this function at the points of a plane $\pi$, as shown in Fig(\ref{plane}). The interpolation is performed for $(x_{1},x_{2},\pi(x_{1},x_{2})$, with $(x_{1},x_{2}) \in [1,18]\times[1,18]$ and $\pi$ is the plane through point $(10.5,10.5,10.5)$ with normal vector equal to $\left(\frac{1}{\sqrt{2}},0,\frac{1}{\sqrt{2}}\right)$. To this end, we utilize the proposed Hermite interpolation as splines $B=B_1 \times B_2 \times B_3$. \textcolor{black}{The other interpolation methods under comparison include Linear spline, cubic spline, Makima (modified Akima cubic Hermite interpolation, from the family of generalized interpolations we implemented b-spline of degree 7 and MOMs of degree 4 and 7. From the classic convolution-based methods we implemented the 6-point Lagrange and 6-point and 8-point cubic polynomial interpolation}.
\textcolor{black}{In order to assess further the behavior of the methods, we also experimented with the density of the 3D grid of support points: we kept the domain constant in each of the $3$ axis and we allowed the step of the grid of the known points to vary from $0.25$ to $1$. As it can be observed from Table \ref{table888}, all interpolation methods (including the 4 tested Hermite polynomials) exhibit reduced error as the step of the grid decreases (equivalently the density of the support points increases). The proposed Hermite polynomial appears to outperform the other methods in this particular example, in terms of RMSE. Furthermore, as expected, the increased spline size ($5 \times 5 \times 5$ vs $3 \times 3 \times 3$) and the higher values of the multiplicity function $\nu$ (2 vs 3) of the proposed Hermite spline result in lower RMSE. Considering the other methods under comparison, the generalized convolution interpolations outperform the convolution-based ones. The differences between the methods are magnified at lower support point density.}

\textcolor{black}{
The behavior of the methods with different grid sizes is visualized in Fig.\ref{plane}(c), where the error with respect to the correct function values as shown in Fig.\ref{plane}(b), is plotted for the proposed Hermite (spline size: 5 and multiplicity function $\nu$ equal to 3), the OMOM of degree 4 and the linear spline, for grid step equal to 0.5 and 1 (upper and lower row respectively). The proposed Hermite interpolant appears to follow the function more closely, which may be attributed to the utilization of the derivatives information}.

\begin{table}[!h]
\begin{center}
\textcolor{black}{
\caption{\textcolor{black}{The achieved RMSE for the proposed Hermite spline polynomial and other interpolation methods, in the case of Example \ref{examp8}}\label{table888}}
\begin{tabular}{|l |cccc|}
\hline
\textbf{Method} &  \multicolumn{4}{c|}{\textbf{Step of the grid}} \\ 
& 0.25  & 0.5 & 0.75 & 1\\ \hline
Hermite $5\times 5 \times 5$, $v(\mathbf{a}) = (2,2,2)$	&2.90e-06 &	0.0021 & 0.0614 &	0.5406\\ \hline
Hermite $5\times 5 \times 5$, $v(\mathbf{a}) = (3,3,3)$	&8.27e-10	& 4.99e-07	&9.21e-05	&0.0049\\ \hline
Hermite $3\times 3 \times 3$, $v(\mathbf{a}) = (2,2,2)$	&6.89e-04	& 0.0355	&0.2690	&1.3867\\ \hline
Hermite $3\times 3 \times 3$, $v(\mathbf{a}) = (3,3,3)$	&  1.34 e-07    & 2.63e-04	& 0.0070 &	0.0908\\ \hline
OMOMS $\deg 4$ 	&2.01e-03 & 0.1549 &	2.3366	&5.7276\\ \hline
OMOMS $\deg 7$ 	& 5.84e-03     &	0.0683 & 1.0745	& 5.5602   \\ \hline
B-spline $\deg 7$ &   5.80e-03	  &	0.0591 &1.0787	&5.4837\\ \hline
Cubic Spline 6p	& 0.0544& 2.2847 &	2.3734	&5.21262\\ \hline
Cubic Spline 8p	& 0.9726 &	0.2756	& 1.5178 &	4.3010\\ \hline
Lagrange 6p	& 0.0261&	1.5501 &	2.5616	&5.1176\\ \hline
Linear spline	&0.5071 &	1.9670	&3.6143	&5.3497\\ \hline
Cubic spline $C^1$	&8.70e-05&	1.0489 &	3.0353	&5.3414\\ \hline
Makima	&1.18e-04 &	1.2054 &	3.1658	&4.9920\\ \hline
Cubic spline $C^2$	&3.17e-06 &	0.4332	& 2.7713 &	5.3440\\ \hline
\end{tabular}
}
\end{center} 
\end{table}
%%%%%%%%%%%%%%%%%%%%%%%%%%%%%%%%%%%%%

\begin{figure}[h!]
%\ContinuedFloat
\centering
\subfigure[]{{\includegraphics[width=0.5\textwidth]{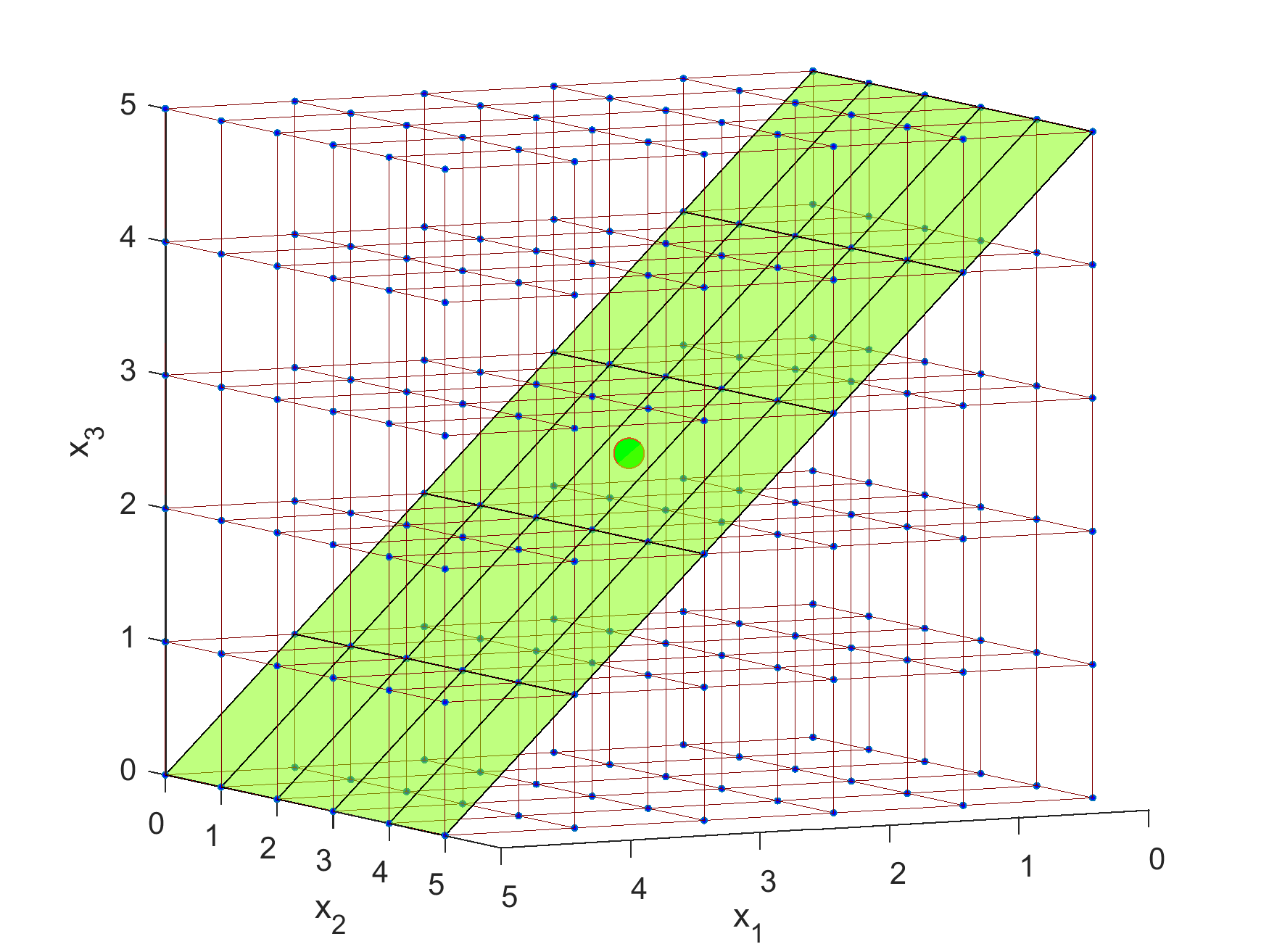}}}\hfill
\subfigure[]{{\includegraphics[width=0.5\textwidth]{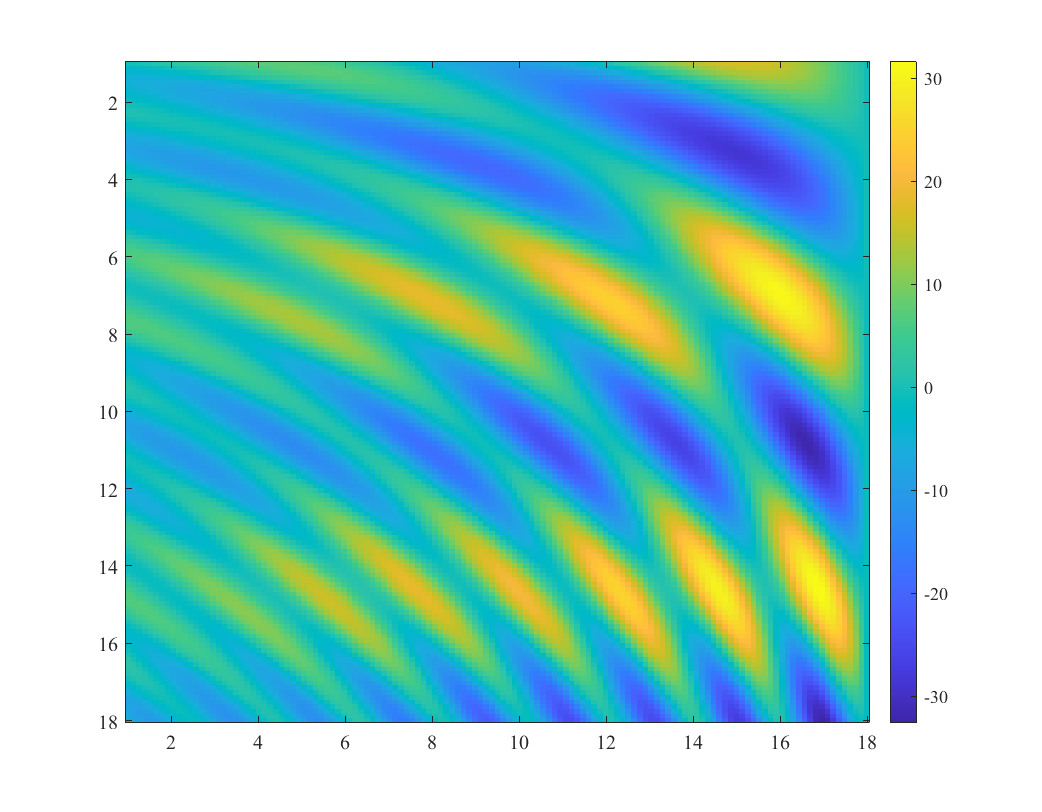}}}\hfill
\subfigure[]{{\includegraphics[width=1\textwidth]{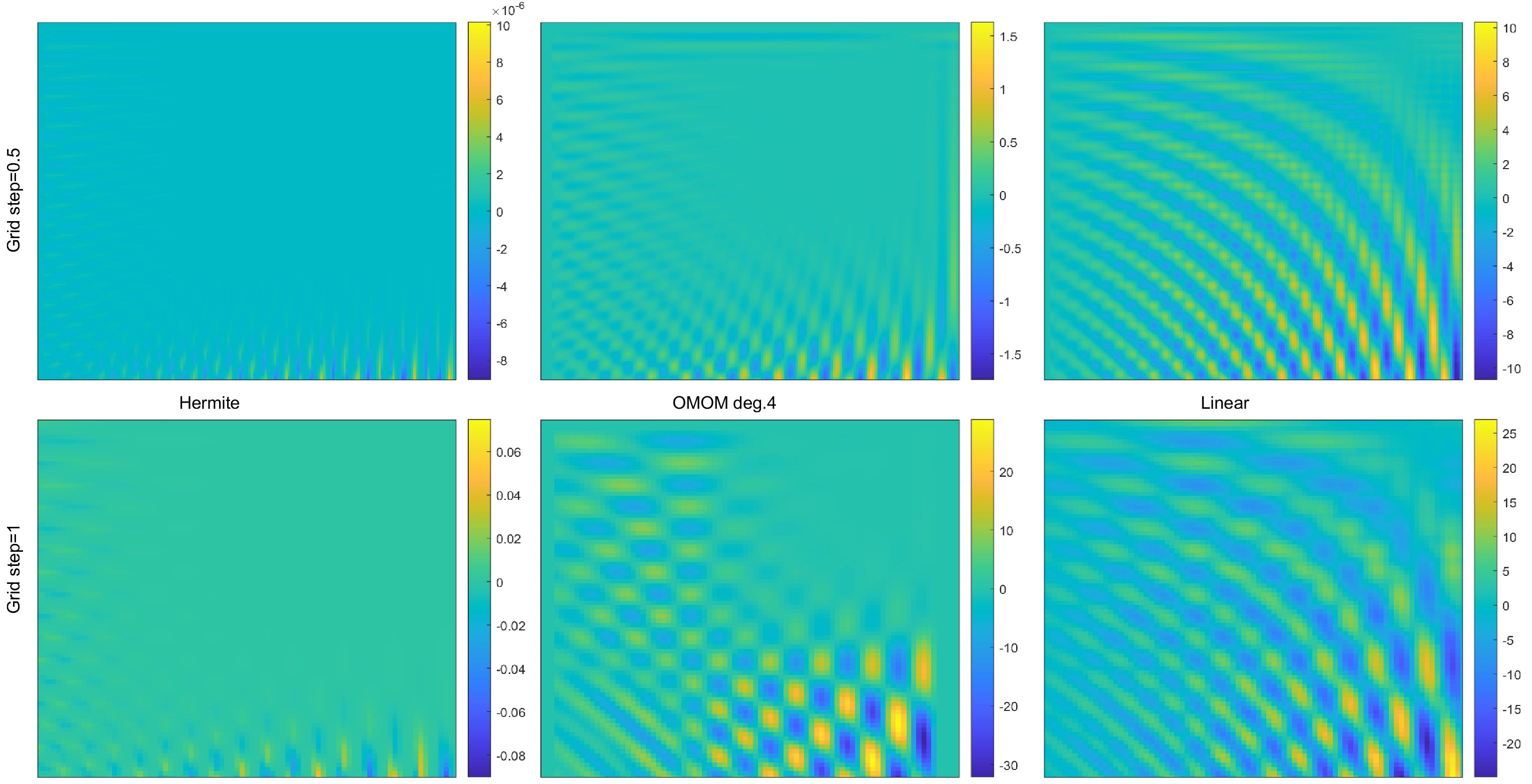}}}
\caption{\textcolor{black}{Spline interpolation along an oblique plane of function $ g(x_1,x_2,x_3) =  x_1 \sin(x_2)+ \frac{x_2 \sin(x_1)}{10}-x_1 \sin\left(\frac{x_2 x_3}{4}\right)$
(a)  the first $6 \times 6\times 6$ support points are displayed, for clarity, with the oblique plane superimposed (b) the actual values of the function $g$ along the oblique plane and (c) the interpolation error (difference from the actual function values) achieved by the proposed Hermite interpolation spline with 5 local support points, multiplicity  $\nu(\mathbf{a})=(3,3,3)$ and grid step 0.5 (upper leftmost image),  Hermite interpolation spline with 5 local support points, multiplicity  $\nu(\mathbf{a})=(3,3,3)$ and grid step 1 (lower-leftmost image),  as well as OMOM deg. 4 with grid step 0.5 (upper center image), OMOM deg. 4 with grid step 1 (lower center image), Linear Spline with grid step 0.5(upper rightmost image), and Linear Spline with grid step 1 (lower rightmost image).}}
\label{plane}
\end{figure}

}
\end{exam}

\newpage
\section{Conclusions}
\label{conc}
\textcolor{black}{In this work we derived a new coordinate Hermite interpolating polynomial,} 
%In this work a new derivation of the coordinate Hermite interpolating polynomial,
given the values of the function and its partial derivatives on $n$-dimensional non-regular grids. The derived formula is algebraically significantly simpler compared to the only alternative formulation in the literature \cite{fractmult}.
Furthermore, we provided the remainder of the interpolation in integral form, as well as cascaded polynomial divisions, in case of interpolating polynomial functions, using the Gr{\"o}bner basis of the ideal of the interpolation. \textcolor{black}{The continuity of the proposed Hermite intepolant polynomials, defined on adjacent $n$-dimensional grids is also proven, establishing spline behavior. The proposed interpolant polynomial can be used for accurate interpolation of data in any number of dimensions, with any degree of partial derivatives at the measured data points. Initial arithmetic results of analytic functions are provided in 2 and 3 dimensions, demonstrating the correctness of the proposed Hermite interpolation. In these examples the proposed Hermite interpolation achieved lower RMSE compared to a number of other interpolation techniques}.

Further work includes the computational implementation of the proposed Hermite polynomial formula utilizing appropriate \textcolor{black}{look-up tables for equidistant points, that will allow fast interpolation of multidimensional signals, as well as a detailed assessment of its accuracy in real-world tasks.}

%% The Appendices part is started with the command \appendix;
%% appendix sections are then done as normal sections
%% \appendix

%% \section{}
%% \label{}

%% References
%%
%% Following citation commands can be used in the body text:
%% Usage of \cite is as follows:
%%   \cite{key}         ==>>  [#]
%%   \cite[chap. 2]{key} ==>> [#, chap. 2]
%%

%% References with BibTeX database:

\bibliographystyle{elsarticle-num}
\bibliography{<your-bib-database>}

\begin{thebibliography}{unsrt}

%% \bibitem must have the following form:
%%   \bibitem{key}...
%%
%%%%%%%%%%%%%%%%%%%%%%%%%%%%%%%%%%%%%
\bibitem{fractmult}{Pandey, K. K., \& Viswanathan, P. (2022). In reference to a self-referential approach towards smooth multivariate approximation. Numerical Algorithms, 91, 251-281. https://doi.org/10.1007/s11075-022-01261-7}

\bibitem{Lorenzt}{Lorentz, R. A. (2000). Multivariate Hermite interpolation by algebraic polynomials: A survey. Journal of Computational and Applied Mathematics, 122(1-2), 167-201. https://doi.org/10.1016/S0377-0427(00)00367-8.}

\bibitem{sauerxu}{Sauer, T., Xu, Y. On multivariate Hermite interpolation. Adv Comput Math 4, 207–259 (1995).}
%%%%%%%%%%%%%%%%%%%%%%%%%%%%%%%%%%%%%%
\bibitem{Birkhoff}{Allasia, G., Cavoretto, R., \& De Rossi, A. (2018). Hermite–Birkhoff interpolation on scattered data on the sphere and other manifolds. Applied Mathematics and Computation, 318, 35-50}.



\bibitem{LR}{Conti, C., Romani, L., \& Unser, M. (2015). Ellipse-preserving Hermite interpolation and subdivision. Journal of Mathematical Analysis and Applications, 426(1), 211-227. https://doi.org/10.1016/j.jmaa.2015.01.017.}

\bibitem{lucia2}{
Romani, L., \& Viscardi, A. (2020). On the refinement matrix mask of interpolating Hermite splines. Applied Mathematics Letters, 109, 106524}



\bibitem{Gasca}{Gasca, M., \& Sauer, T. (2000). On Bivariate Hermite Interpolation with Minimal Degree Polynomials. SIAM Journal on Numerical Analysis, 37(3), 772–798.} http://www.jstor.org/stable/2587315

\bibitem{Gevorgian}{Gevorgian, H.V., Hakopian, H.A. \& Sahakian, A.A. On the bivariate hermite interpolation problem. Constr. Approx 11, 23–35 (1995).} https://doi.org/10.1007/BF01294336

\bibitem{Salzerbivhyperosculatory}{H.E. Salzer, Formulas for bivariate hyperosculatory interpolation, Math. Comput. 25 (1971) 119–133}

\bibitem{GascaSauerhistory}{M. Gasca, T. Sauer On the history of multivariate polynomial interpolation, Journal of Computational and Applied Mathematics 122 (2000) 23–35 33}

\bibitem{unser1} {Unser, M., Aldroubi, A., \& Eden, M. (1993). B-spline signal processing. I. Theory. IEEE transactions on signal processing, 41(2), 821-833.}

\bibitem{unser2}{Unser, M., Aldroubi, A., \& Eden, M. (1993). B-spline signal processing. II. Efficiency design and applications. IEEE transactions on signal processing, 41(2), 834-848.}

\bibitem{omom1}{Blu, T., Th\'evenaz, P., \& Unser, M. (2001). MOMS: Maximal-order interpolation of minimal support. IEEE Transactions on Image Processing, 10(7), 1069-1080.}

\bibitem{interp_revisit}{Th\'evenaz, P., Blu, T., \& Unser, M. (2000). Interpolation revisited [medical images application]. IEEE Transactions on medical imaging, 19(7), 739-758.}

\bibitem{SPITZ}{Spitzbart, A. (1960). A Generalization of Hermite’s Interpolation Formula. The American Mathematical Monthly, 67(1), 42–46. https://doi.org/10.2307/2308924}

\bibitem{DK12} {Delibasis, K., Kechriniotis, A., \& Assimakis, N. (2012). New closed formula for the univariate Hermite interpolating polynomial of total degree and its application in medical image slice interpolation. IEEE Transactions on Signal Processing, 60(12), 6294-6304.}

\bibitem{GHITV}{Chawla, M., \& Jayarajan, N. (1974). A generalization of Hermite's interpolation formula in two variables. Journal of the Australian Mathematical Society, 18(4), 402-410.}

\bibitem{gascasa}{Gasca, M., \& Sauer, T. (2000). Polynomial interpolation in several variables. Advances in Computational Mathematics, 12, 377-410.}

\bibitem{Sauer2006}{
Sauer, T. (2006). Polynomial interpolation in several variables: Lattices, differences, and ideals. In M. Buhmann, W. Hausmann, K. Jetter, W. Schaback, \& J. Stöckler (Eds.), Multivariate Approximation and Interpolation (pp. 189-228). Elsevier.}






\bibitem{DK14} {Delibasis, K. K., \& Kechriniotis, A. (2014). A new formula for bivariate Hermite interpolation on variable step grids and its application to image interpolation. IEEE Transactions on Image Processing, 23(7), 2892-2904. https://doi.org/10.1109/TIP.2014.2322441.}

\bibitem{BGHI} {Ahlin, A. C. (1964). A bivariate generalization of Hermite's interpolation formula. Mathematics of Computation, 18(86), 264-273.}


\bibitem{Davis1975}{
Davis, P. J. (1975). Interpolation and approximation. Dover Books on Advanced Mathematics. Dover Publications.}

\bibitem{Vladimirov2007}{
Vladimirov, V. S. (2007). Methods of the theory of functions of many complex variables. MIT Press. (Original work published 1966) Dover reprint.}

\bibitem{Fisher1999}{
Fisher, S. D. (1999). Complex variables. Wadsworth \& Brooks. (Original work published 1990)}

\bibitem{deBoor2005}{
de Boor, C. (2005). Ideal interpolation. In C. K. Chui, M. Neamtu, \& L. L. Schumaker (Eds.), Approximation Theory XI, Gaitlinburg 2004 (pp. 59-91). Nashboro Press.}

\bibitem{ideal_Shekhtman}{Shekhtman, B. (2009). Ideal interpolation: translations to and from algebraic geometry. In Approximate commutative algebra (pp. 163-192). Vienna: Springer Vienna.}

\bibitem{Sauer2018}{
Sauer, T. (2018). Prony's method in several variables: symbolic solutions by universal interpolation. Journal of Symbolic Computation, 84, 95-112. arXiv:1603.03944.}


%\bibitem{ideal_Bazan2023}{Bazan, E. R., \& Hubert, E. (2023). Symmetry in multivariate ideal interpolation. Journal of Symbolic Computation, 115, 174-200.}

%\bibitem{Barthelmann2000}{Barthelmann, V., Novak, E., \& Ritter, K. (2000). High dimensional polynomial interpolation on sparse grids. Advances in Computational Mathematics, 12, 273-288. Special Issue.}



%\bibitem{Carnicer2006}{Carnicer, J., Gasca, M., \& Sauer, T. (2006). Interpolation lattices in several variables. Numerische Mathematik, 102, 559-581.}

\bibitem{practspline}{de Boor, C. (1978). A Practical Guide to Splines. Springer-Verlag.}













%%%%%%%%%%%%%%%%%%%%%%%%%%%%%%%%%%%%%%%%%%%%%%


\end{thebibliography}

%% Authors are advised to use a BibTeX database file for their reference list.
%% The provided style file elsarticle-num.bst formats references in the required Procedia style

%% For references without a BibTeX database:

\end{document}